\documentclass[12pt]{article}
\usepackage{amssymb, amscd, amsmath, amsthm, amsfonts, ae, color}
\usepackage[all,2cell,dvips]{xy} 

\title{Realizations of $BC_r$-graded intersection matrix algebras with grading subalgebras of type $B_r$, $r \geq 3$}
\author{\emph{In memory of Professor Peter Slodowy}}
\date{S. Bhargava, \ Y. Gao\footnote{Funding from the National Sciences and Engineering Research Council of Canada is gratefully acknowledged.}}
\begin{document}


\definecolor{darkslateblue}{rgb}{0.284,0.240,0.545}
\definecolor{brightblue}{rgb}{0.000,0.400,1.000}
\definecolor{midnightblue}{rgb}{0.200,0.200,0.600}
\definecolor{lightblue}{rgb}{0.680,0.848,0.900}
\definecolor{lightskyblue}{rgb}{0.530,0.808,0.980}
\definecolor{powderblue}{rgb}{0.880,0.880,0.950}

\definecolor{redbrown}{rgb}{0.800,0.400,0.200}

\definecolor{green}{rgb}{0,.5,0}
\definecolor{greeen}{rgb}{.1,.7,.1}
\definecolor{greeeen}{rgb}{.1,.5,.1}
\definecolor{darkgreen}{rgb}{0.000,0.392,0.000}
\definecolor{forestgreen}{rgb}{.133333,.545098,.133333}
\definecolor{grassgreen}{rgb}{0.200,0.600,0.000}
\definecolor{springgreen}{rgb}{0.000,1.000,0.498}

\definecolor{lightgrey}{rgb}{0.828,0.828,0.828}

\definecolor{royalpurple}{rgb}{0.400,0.000,0.600}

\definecolor{redorange}{rgb}{0.800,0.200,0.000}
\definecolor{darkorange}{rgb}{1.000,0.550,0.000}

\definecolor{darkred}{rgb}{0.545098,0,0}
\definecolor{tomato}{rgb}{1.000,0.390,0.280}
\definecolor{firebrick}{rgb}{0.698,0.132,0.132}
\definecolor{indianred}{rgb}{0.804,0.360,0.360}

\definecolor{darkviolet}{rgb}{.580392,0,.827451}

\definecolor{softyellow}{rgb}{1.000,1.000,0.500}

\newtheorem{thm}{Theorem}
\newtheorem*{thm*}{Theorem}
\newtheorem{defn}{Definition}
\newtheorem{propn}{Proposition}
\newtheorem*{propn*}{Proposition}
\newtheorem{lemma}{Lemma}
\newtheorem*{lemma*}{Lemma}
\newtheorem{cor}{Corollary}
\newtheorem*{cor*}{Corollary}
\newtheorem{eg}{Example}
\newtheorem{claim}{Claim}

\newcommand{\F}{\mathbb{F}}
\newcommand{\C}{\mathbb{C}}
\newcommand{\R}{\mathbb{R}}
\newcommand{\bbC}{\mathbb{C}}
\newcommand{\bbF}{\mathbb{F}}
\newcommand{\bbQ}{\mathbb{Q}}
\newcommand{\bbR}{\mathbb{R}}
\newcommand{\bbZ}{\mathbb{Z}}

\newcommand{\calA}{\mathcal{A}}
\newcommand{\calB}{\mathcal{B}}
\newcommand{\calE}{\mathcal{E}}
\newcommand{\calG}{\mathcal{G}}
\newcommand{\calH}{\mathcal{H}}
\newcommand{\calL}{\mathcal{L}}
\newcommand{\calM}{\mathcal{M}}
\newcommand{\calS}{\mathcal{S}}

\newcommand{\frake}{\mathfrak{e}}
\newcommand{\frakg}{\mathfrak{g}}
\newcommand{\mfa}{\mathfrak{a}}
\newcommand{\mfb}{\mathfrak{b}}
\newcommand{\mfc}{\mathfrak{c}}
\newcommand{\mfg}{\mathfrak{g}}
\newcommand{\mfh}{\mathfrak{h}}
\newcommand{\mfi}{\mathfrak{i}}
\newcommand{\mfr}{\mathfrak{r}}
\newcommand{\mfz}{\mathfrak{z}}
\newcommand{\mfA}{\mathfrak{A}}
\newcommand{\mfB}{\mathfrak{B}}

\newcommand{\lcb}{\{} 
\newcommand{\rcb}{\}} 

\newcommand{\fd}{\texttt{finite-dimensional }}
\newcommand{\dm}{\text{dim}}  
\newcommand{\End}{\text{End}}
\newcommand{\EndFV}{\text{End}_{\mathbb{F}}(V)}
\newcommand{\ns}{\text{nullspace}} 
\newcommand{\spn}{\text{span}} 

\newcommand{\ad}{\text{ad} \,}
\newcommand{\gim}{\mathfrak{g}\mathfrak{i}\mathfrak{m}}
\newcommand{\gimAd}{\mathfrak{g}\mathfrak{i}\mathfrak{m}\!\left( A^{[d]} \right)}
\newcommand{\gimSLE}{\mathfrak{g}\mathfrak{i}\mathfrak{m}\big( A^{[\mathcal{S},\mathcal{L},\mathcal{E}]}  \big)}
\newcommand{\im}{\mathfrak{i}\mathfrak{m}}
\newcommand{\imAd}{\mathfrak{i}\mathfrak{m}\!\left( A^{[d]} \right)}

\newcommand{\Mtn}{\text{M}_{2n}} 
\newcommand{\MtnR}{\text{M}_{2n}(R)} 
\newcommand{\gltnR}{\text{gl}_{2n}(R)} 
\newcommand{\sotnpoaecc}{\text{so}_{2r+1}\left(\mathfrak{a},\eta,C,\chi\right)}
\newcommand{\sotrpoaecc}{\text{so}_{2r+1}\left(\mathfrak{a},\eta,C,\chi\right)}
\newcommand{\sotrpobecc}{\text{so}_{2r+1}\left(\mathfrak{b},\eta,C,\chi\right)}
\newcommand{\sptnC}{\text{sp}_{2n}(\mathbb{C})}
\newcommand{\sptnR}{\text{sp}_{2n}(R, -)}
\newcommand{\spe}{\text{sp}_8 \big( FG(x,y,z), - \big) }

\newcommand{\uinv}{u^{-1}}
\newcommand{\vinv}{v^{-1}}
\newcommand{\winv}{w^{-1}}
\newcommand{\xinv}{x^{-1}}
\newcommand{\yinv}{y^{-1}}
\newcommand{\zinv}{z^{-1}}

\newcommand{\Tone}{t_{\mu_1,k_1}} 
\newcommand{\Tonem}{t_{\mu_1,k_1}^{m_1}}
\newcommand{\Toneinv}{t_{\mu_1,k_1}^{-1}}
\newcommand{\Tmanyl}{t_{\mu_1,k_1}^{m_1} \cdots t_{\mu_l,k_l}^{m_l}}
\newcommand{\Tlpo}{t_{\mu_{l+1},k_{l+1}}}
\newcommand{\Tlpom}{t_{\mu_{l+1},k_{l+1}}^{m_{l+1}}}
\newcommand{\Tmanylpo}{t_{\mu_1,k_1}^{m_1} \cdots t_{\mu_l,k_l}^{m_l} t_{\mu_{l+1},k_{l+1}}^{m_{l+1}}}
\newcommand{\xmuk}{x_{\mu_1,k_1}}
\newcommand{\xmukinv}{x_{\mu_1,k_1}^{-1}}
\newcommand{\ymuk}{y_{\mu_1,k_1}}
\newcommand{\ymukinv}{y_{\mu_1,k_1}^{-1}}
\newcommand{\zmuk}{z_{\mu_1,k_1}}
\newcommand{\zmukinv}{z_{\mu_1,k_1}^{-1}}

\newcommand{\Alps}{A^{[l+s]}}

\newcommand{\Eijr}{E_{ij}(r)}
\newcommand{\Enjnir}{E_{n+j,n+i}(\bar{r})}
\newcommand{\Einir}{E_{i,n+i}(r)}
\newcommand{\Einirbar}{E_{i,n+i}(\bar{r})}
\newcommand{\Einjr}{E_{i,n+j}(r)}
\newcommand{\Ejnirbar}{E_{j,n+i}(\bar{r})}
\newcommand{\Eniir}{E_{n+i,i}(r)}
\newcommand{\Eniirbar}{E_{n+i,i}(\bar{r})}
\newcommand{\Enijr}{E_{n+i,j}(r)}
\newcommand{\Enjirbar}{E_{n+j,i}(\bar{r})}
\newcommand{\Eiir}{E_{ii}(r)}
\newcommand{\Eninirbar}{E_{n+i,n+i}(\bar{r})}

\newcommand{\eii}{e_{ii}}
\newcommand{\ejj}{e_{jj}}
\newcommand{\enini}{e_{n+i,n+i}}
\newcommand{\enjnj}{e_{n+j,n+j}}

\newcommand{\Cg}{\mathfrak{C}\left(\mathfrak{g}\right)}

\newcommand{\Ginv}{G^{-1}}
\newcommand{\Jr}{J_r}
\newcommand{\Mtrpo}{M_{2r+1}}
\newcommand{\ncLpxy}{\mathbb{C}\left[x^{\pm 1}, y^{\pm 1} \right]_{\text{n.c.}} } 
\newcommand{\sotrpo}{\text{so}_{2r+1}}
\newcommand{\sotrpoC}{\text{so}_{2r+1}\left(\mathbb{C}\right)}
\newcommand{\sotrpobe}{\text{so}_{2r+1}\left(\mathfrak{b}, \eta \right)}

\newcommand{\tinyp}{$\phantom{\tiny p}$}
\newcommand{\blobbb}{\rule[-.2ex]{1ex}{1ex}}
\newcommand{\Eul}{E^{\hbox{
  \begin{picture}(10,10)
    \put(0,0){\framebox(8,8)[0,0]{}}
    \linethickness{4pt}  \put(0,6){\line(1,0){4}}
  \end{picture}}}}
\newcommand{\Eur}{E^{\hbox{
  \begin{picture}(10,10)
    \put(0,0){\framebox(8,8)[0,0]{}}
    \linethickness{4pt}  \put(4,6){\line(1,0){4}}
  \end{picture}}}}
\newcommand{\Ebl}{E^{\hbox{
  \begin{picture}(10,10)
    \put(0,0){\framebox(8,8)[0,0]{}}
    \linethickness{4pt}  \put(0,2){\line(1,0){4}}
  \end{picture}}}}

\newcommand{\ulbox}{\hbox{
  \begin{picture}(10,10)
    \put(0,0){\framebox(8,8)[0,0]{}}
    \linethickness{4pt}  \put(0,6){\line(1,0){4}}
  \end{picture}}}
\newcommand{\urbox}[2]{\hbox{
  \begin{picture}(10,10)
    \put(0,0){\framebox(8,8)[0,0]{}}
    \linethickness{4pt}  \put(4,6){\line(1,0){4}}
  \end{picture}}}
\newcommand{\blbox}{\hbox{
  \begin{picture}(10,10)
    \put(0,0){\framebox(8,8)[0,0]{}}
    \linethickness{4pt}  \put(0,2){\line(1,0){4}}
  \end{picture}}}
\newcommand{\hortbox}{\hbox{
  \begin{picture}(10,10)
    \put(0,0){\framebox(8,8)[0,0]{}}
    \put(0,4){\line(1,0){4}}
  \end{picture}}}
\newcommand{\vertbox}{\hbox{
  \begin{picture}(10,10)
    \put(0,0){\framebox(8,8)[0,0]{}}
    \put(4,4){\line(0,1){4}}
  \end{picture}}}

\newcommand{\eppo}{\epsilon_{p+1}}
\newcommand{\EndaC}{\text{End}_{\mathfrak{a}}(C)}
\newcommand{\EndanpC}{\text{End}_{\mathfrak{a}} \left( \mathfrak{a}^n \oplus C  \right)}
\newcommand{\keiminus}{k_{e_i}^-}
\newcommand{\keiplus}{k_{e_i}^+}
\newcommand{\kli}{k_{l_i}}
\newcommand{\mfanpC}{\mathfrak{a}^n \oplus C}
\newcommand{\Mna}{\text{M}_n (\mathfrak{a})}
\newcommand{\Nalpha}{N_{\alpha}}
\newcommand{\Nomega}{N_{\omega}}
\newcommand{\Ntheta}{N_{\theta}}
\newcommand{\ovla}{\overline{a}}
\newcommand{\ovlb}{\overline{b}}
\newcommand{\ovlab}{\overline{a}\,\overline{b}}
\newcommand{\ovlba}{\overline{b}\,\overline{a}}
\newcommand{\ovlp}{\overline{p}}
\newcommand{\ovlq}{\overline{q}}
\newcommand{\ovlu}{\overline{u}}
\newcommand{\ovlv}{\overline{v}}
\newcommand{\ovlw}{\overline{w}}
\newcommand{\ovluv}{\overline{u}\ \overline{v}}
\newcommand{\ovlvu}{\overline{v}\ \overline{u}}
\newcommand{\ovluw}{\overline{u}\ \overline{w}}
\newcommand{\ovlwu}{\overline{w}\ \overline{u}}
\newcommand{\ovlvw}{\overline{v}\ \overline{w}}
\newcommand{\ovlwv}{\overline{w}\ \overline{v}}
\newcommand{\ovluinv}{\overline{u^{-1}}}
\newcommand{\ovlvinv}{\overline{v^{-1}}}
\newcommand{\ovlwinv}{\overline{w^{-1}}}
\newcommand{\ovlx}{\overline{x}}
\newcommand{\ovly}{\overline{y}}
\newcommand{\ovlz}{\overline{z}}
\newcommand{\ovlxy}{\overline{x}\,\overline{y}}
\newcommand{\ovlyx}{\overline{y}\,\overline{x}}
\newcommand{\ovlxinv}{\overline{x^{-1}}}
\newcommand{\ovlyinv}{\overline{y^{-1}}}
\newcommand{\TIl}{\text{Type I}\,\mathit{l}}
\newcommand{\TIs}{\text{Type I}\,\mathit{s}}
\newcommand{\TIIl}{\text{Type II}\,\mathit{l}}
\newcommand{\TIIs}{\text{Type II}\,\mathit{s}}
\newcommand{\TIIeplus}{\text{Type II}\,\mathit{e}^+}
\newcommand{\TIIeminus}{\text{Type II}\,\mathit{e}^-}
\newcommand{\usc}{\underline{c}}
\newcommand{\usz}{\underline{z}}

\newcommand{\Eone}{E_{1,2}(1) + E_{8,9}(-1)}
\newcommand{\Etwo}{E_{2,3}(1) + E_{7,8}(-1)}
\newcommand{\Ethree}{E_{3,4}(1) + E_{6,7}(-1)}
\newcommand{\Efour}{E_{4,5}(\sqrt{2}) + E_{5,6}(-\sqrt{2})}
\newcommand{\Efive}{E_{5,1}(p) + E_{9,5}(-\overline{p}) + E_{9,10}\big( \sqrt{2} \, \chi_{m} \big) + E_{10,1}(\sqrt{2} \, m)  }
\newcommand{\EEfive}{E_{9,10}\big( \sqrt{2} \, \chi_{m} \big) + E_{10,1}(\sqrt{2} \, m)  }

\newcommand{\Fone}{E_{2,1}(1) + E_{9,8}(-1)}
\newcommand{\Ftwo}{E_{3,2}(1) + E_{8,7}(-1)}
\newcommand{\Fthree}{E_{4,3}(1) + E_{7,6}(-1)}
\newcommand{\Ffour}{E_{5,4}(\sqrt{2}) + E_{6,5}(-\sqrt{2})}
\newcommand{\Ffive}{E_{1,5}(q) + E_{5,9}(-\overline{q}) + E_{1,10}\big( \sqrt{2} \, \chi_{n} \big) + E_{10,9}(\sqrt{2} \, n)  }
\newcommand{\FFfive}{E_{1,10}\big( \sqrt{2} \, \chi_{n} \big) + E_{10,9}(\sqrt{2} \, n)  }

\newcommand{\Hone}{E_{1,1}(1) + E_{9,9}(-1) + E_{2,2}(-1) + E_{8,8}(1)}
\newcommand{\Htwo}{E_{2,2}(1) + E_{8,8}(-1) + E_{3,3}(-1) + E_{7,7}(1)}
\newcommand{\Hthree}{E_{3,3}(1) + E_{7,7}(-1) + E_{4,4}(-1) + E_{6,6}(1)}
\newcommand{\Hfour}{E_{4,4}(2) + E_{6,6}(-2)}
\newcommand{\HHfive}{E_{1,1}\big( - 2 \, \chi(n,m) \big) + E_{9,9}\big( 2 \, \chi(m,n) \big) + E_{10,10}\big( 2 \, m \, \chi_n - 2 \, n \, \chi_m \big)} 
\newcommand{\HHHfive}{E_{1,1}( - 2 ) + E_{9,9}( 2 ) + E_{10,10}\big( 2 \, m \, \chi_n - 2 \, n \, \chi_m \big)} 

\newcommand{\te}{\tilde{e}}
\newcommand{\tf}{\tilde{f}}
\newcommand{\tilh}{\tilde{h}}
\newcommand{\tL}{\widetilde{L}}
\newcommand{\tU}{\widetilde{U}}

\newcommand{\UuL}{U_{\underline{L}}}

\newcommand{\eimej}[2]{\epsilon_{#1} - \epsilon_{#2} }
\newcommand{\eipej}[2]{\epsilon_{#1} + \epsilon_{#2} }
\newcommand{\meimej}[2]{-\epsilon_{#1} - \epsilon_{#2} }

\newcommand{\Eulquad}[3]{E_{#1,#2}(#3) + E_{2r+2-#2,2r+2-#1}(-#3)} 
\newcommand{\EULquad}[6]{E_{#1,#2}(#3) + E_{2r+2-#4,2r+2-#5}(#6)}
\newcommand{\Eurquad}[3]{E_{#1,2r+2-#2}(#3) + E_{#2,2r+2-#1}(-#3)} 
\newcommand{\EURquad}[6]{E_{#1,2r+2-#2}(#3) + E_{#4,2r+2-#5}(#6)}
\newcommand{\Eblquad}[3]{E_{2r+2-#1,#2}(#3) + E_{2r+2-#2,#1}(-#3)}
\newcommand{\EBLquad}[6]{E_{2r+2-#1,#2}(#3) + E_{2r+2-#4,#5}(#6)}
\newcommand{\Eplusr}[2]{E_{#1,r+1}(#2) + E_{r+1,2r+2-#1}(-#2)} 
\newcommand{\EPLUSR}[3]{E_{#1,r+1}(#2) + E_{r+1,2r+2-#1}(#3)} 
\newcommand{\Eminusr}[2]{E_{r+1,#1}(#2) + E_{2r+2-#1,r+1}(-#2)} 
\newcommand{\EMINUSR}[3]{E_{r+1,#1}(#2) + E_{2r+2-#1,r+1}(#3)} 
\newcommand{\Ediag}[2]{E_{#1,#1}(#2) + E_{2r+2-#1,2r+2-#1}(-#2)} 
\newcommand{\EDIAG}[3]{E_{#1,#1}(#2) + E_{2r+2-#1,2r+2-#1}(#3)}
\newcommand{\EDIAGbrac}[3]{E_{#1,#1}(#2) + E_{2r+2-(#1),2r+2-(#1)}(#3)}

\newcommand{\Eeimej}[3]{E_{\epsilon_{#1} - \epsilon_{#2}} (#3) }
\newcommand{\Eeipej}[3]{E_{\epsilon_{#1} + \epsilon_{#2}} (#3) }
\newcommand{\Emeimej}[3]{E_{-\epsilon_{#1} - \epsilon_{#2}} (#3) }
\newcommand{\Eei}[2]{E_{\epsilon_{#1}} (#2) }
\newcommand{\Emei}[2]{E_{-\epsilon_{#1}} (#2) }
\newcommand{\Hi}[2]{H_{#1} (#2) }

\newcommand{\Meimej}[2]{M_{\epsilon_{#1} - \epsilon_{#2}} }
\newcommand{\Meipej}[2]{M_{\epsilon_{#1} + \epsilon_{#2}} }
\newcommand{\Mmeimej}[2]{M_{\epsilon_{#1} - \epsilon_{#2}} }
\newcommand{\mk}[2]{m_{#1,#2}}
\newcommand{\xmi}[1]{x^{m_{#1,x}}}
\newcommand{\ymi}[1]{y^{m_{#1,y}}}

\newcommand{\ovl}[1]{\overline{#1}}
\newcommand{\teeimej}[2]{\tilde{e}_{\epsilon_{#1} - \epsilon_{#2}} }
\newcommand{\tfeimej}[2]{\tilde{f}_{\epsilon_{#1} - \epsilon_{#2}} }
\newcommand{\theimej}[2]{\tilde{h}_{\epsilon_{#1} - \epsilon_{#2}} }
\newcommand{\teeipej}[2]{\tilde{e}_{\epsilon_{#1} + \epsilon_{#2}} }
\newcommand{\tfeipej}[2]{\tilde{f}_{\epsilon_{#1} + \epsilon_{#2}} }
\newcommand{\theipej}[2]{\tilde{h}_{\epsilon_{#1} + \epsilon_{#2}} }
\newcommand{\temeimej}[2]{\tilde{e}_{-\epsilon_{#1} - \epsilon_{#2}} }
\newcommand{\tfmeimej}[2]{\tilde{f}_{-\epsilon_{#1} - \epsilon_{#2}} }
\newcommand{\thmeimej}[2]{\tilde{h}_{-\epsilon_{#1} - \epsilon_{#2}} }

\newcommand{\EUL}[2]{E^{\hbox{
  \begin{picture}(10,10)
    \put(0,0){\framebox(8,8)[0,0]{}}
    \linethickness{4pt}  \put(0,6){\line(1,0){4}}
  \end{picture}}}_{#1,#2}}

\newcommand{\EUR}[2]{E^{\hbox{
  \begin{picture}(10,10)
    \put(0,0){\framebox(8,8)[0,0]{}}
    \linethickness{4pt}  \put(4,6){\line(1,0){4}}
  \end{picture}}}_{#1,2r+2-#2}}

\newcommand{\EBL}[2]{E^{\hbox{
  \begin{picture}(10,10)
    \put(0,0){\framebox(8,8)[0,0]{}}
    \linethickness{4pt}  \put(0,2){\line(1,0){4}}
  \end{picture}}}_{2r+2-#1,#2}}

\newcommand{\Ehort}[1]{E^{\hbox{
  \begin{picture}(10,10)
    \put(0,0){\framebox(8,8)[0,0]{}}
    \put(0,4){\line(1,0){4}}
  \end{picture}}}_{r+1,#1}}

\newcommand{\Evert}[1]{E^{\hbox{
  \begin{picture}(10,10)
    \put(0,0){\framebox(8,8)[0,0]{}}
    \put(4,4){\line(0,1){4}}
  \end{picture}}}_{#1,r+1}}

\maketitle

\abstract{
We study intersection matrix algebras $\mathfrak{i}\mathfrak{m}\!\left( A^{[d]} \right)$ that arise from affinizing a Cartan matrix $A$ of type $B_r$ with $d$ arbitrary long roots in the root system $\Delta_{B_r}$, where $r \geq 3$.  We show that $\mathfrak{i}\mathfrak{m}\!\left( A^{[d]} \right)$ is isomorphic to the universal covering algebra of $\text{so}_{2r+1}\left(\mathfrak{a},\eta,C,\chi\right)$, where $\mathfrak{a}$ is an associative algebra with involution $\eta$, and $C$ is an $\mathfrak{a}$-module with hermitian form $\chi$.  We provide a description of all four of the components $\mfa$, $\eta$, $C$, and $\chi$.
}

\vskip10pt

\noindent 2000 Mathematics Subject Classification: Primary 17B65, 17B70, 17B05 \qquad Secondary 17B67, 16W10

\section{Introduction}

In the early to mid-1980s, Peter Slodowy discovered that matrices like 
	\begin{displaymath}
		M =
			\left[
				\begin{array}{rrrr}
					2 & -1 & 0 & 1\\
					-1 & 2 & -1 & 1\\
					0 & -2 & 2 & -2 \\
					1 & 1 & -1 & 2
				\end{array}
			\right]
	\end{displaymath}
were encoding the intersection form on the second homology group of Milnor fibres for germs of holomorphic maps with an isolated singularity at the origin [S1], [S2].  These matrices were like the generalized Cartan matrices of Kac-Moody theory in that they had integer entries, $2$'s along the diagonal, and $M_{ij}$ was negative if and only if $M_{ji}$ was negative.  What was new, however, was the presence of positive entries off the diagonal.  Slodowy called such matrices generalized intersection matrices:

	\begin{defn}[{[S1]}]\emph{
		An $n \times n$ integer-valued matrix $M$ is called a \emph{generalized intersection matrix} ($\mathfrak{g}\mathfrak{i}\mathfrak{m}$) if 
			\begin{itemize}
				\item[]
					$M_{ii} =2$,
				\item[]
					$M_{ij} < 0$ iff $M_{ji} < 0$, and 
				\item[]
					$M_{ij} >0$ iff $M_{ji} > 0$,
			\end{itemize}
		for $1 \leq i,j \leq n$ with $i \neq j$.
	}\end{defn}

Slodowy used these matrices to define a class of Lie algebras that encompassed all the Kac-Moody Lie algebras:

\begin{defn}[{[S1],[BrM]}]\label{defn:gim algebra}\emph{
	Given an $n \times n$ generalized intersection matrix $M = \left( M_{ij} \right)$, define a Lie algebra over $\bbC$, called a \emph{generalized intersection matrix ($\mathfrak{g}\mathfrak{i}\mathfrak{m}$) algebra} and denoted by $\gim(M)$, with:
		\begin{itemize}
			\item[]
				generators:  \qquad $e_1, \ldots, e_n$, $f_1, \ldots, f_n$, $h_1, \ldots h_n$
			\item[]
				relations:
					\begin{itemize}
						\item[(R1)]
							for $1 \leq i,j \leq n$,
								\begin{itemize}
									\item[]
										$[h_i, e_j] = M_{ij} e_j$
									\item[]
										$[h_i, f_j] = -M_{ij} f_j$
									\item[]
										$[e_i, f_i] = h_i$
								\end{itemize}  
						\item[(R2)]
							for $M_{ij} \leq 0$,
								\begin{itemize}
									\item[]
										$[e_i,f_j] = 0 = [f_i,e_j]$
									\item[]
										$(\ad e_i)^{-M_{ij} +1} \ e_j = 0 = (\ad f_i)^{-M_{ij} +1} \ f_j$
								\end{itemize}
						\item[(R3)]
							for $M_{ij} >0$, $i \neq j$
								\begin{itemize}
									\item[]
										$[e_i,e_j] = 0 = [f_i,f_j]$
									\item[]
										$(\ad e_i)^{M_{ij} +1} \ f_j = 0 = (\ad f_i)^{M_{ij} +1} \ e_j$
								\end{itemize}	
					\end{itemize}
				\end{itemize}
	}\end{defn}

If the $M$ that we begin with is a generalized Cartan matrix, then the $3n$ generators and the first two groups of axioms, (R1) and (R2), provide a presentation of the Kac-Moody Lie algebras [GbK], [C], [K].

Slodowy and, later, Berman showed that the $\mathfrak{g}\mathfrak{i}\mathfrak{m}$ algebras are also isomorphic to fixed point subalgebras of involutions on larger Kac-Moody algebras [S1], [Br].  So, in their words, the $\mathfrak{g}\mathfrak{i}\mathfrak{m}$ algebras lie both ``beyond and inside'' Kac-Moody algebras.  

Further progress came in the 1990s as a byproduct of the work of Berman and Moody, Benkart and Zelmanov, and Neher on the classification of root-graded Lie algebras [BrM], [BnZ], [N].  Their work revealed that some families of intersection matrix ($\mathfrak{i}\mathfrak{m}$) algebras, which are quotient algebras of $\mathfrak{g}\mathfrak{i}\mathfrak{m}$ algebras, were universal covering algebras of well understood Lie algebras.  For instance the $\mathfrak{i}\mathfrak{m}$ algebras that arise from multiply affinizing a Cartan matrix of type $A_r$, with $r \geq 3$, are the universal covering algebras of $\text{sl}(\mfa)$, where $\mfa$ is the associative algebra of noncommuting Laurent polynomials in several variables (the number of indeterminates depends on how many times the original Cartan matrix  is affinized).  A handful of other researchers also began engaging these new algebras.  For example, Eswara Rao, Moody, and Yokonuma used vertex operator representations to show that $\mathfrak{i}\mathfrak{m}$ algebras were nontrivial [EMY].  Gao examined compact forms of $\mathfrak{i}\mathfrak{m}$ algebras arising from conjugations over the complex field [G]. Peng found relations between $\im$ algebras and the representations of tilted algebras via Ringel-Hall algebras [P]. Berman, Jurisich, and Tan showed that the presentation of $\mathfrak{g}\mathfrak{i}\mathfrak{m}$ algebras could be put into a broader framework that incorporated Borcherds algebras [BrJT].  

The chief objective of this paper is to continue advancing our understanding of $\mathfrak{g}\mathfrak{i}\mathfrak{m}$ and $\mathfrak{i}\mathfrak{m}$ algebras.  We construct a generalized intersection matrix $A^{[d]}$ by adjoining $d$ long roots to a base of a root system of type $B_r$, where $r \geq 3$.  This is exactly the analogue of the affinization process in which a single root is adjoined to a Cartan matrix of a finite-dimensional Lie algebra to arrive at a generalized Cartan matrix and, eventually, an affine Kac-Moody algebra. The matrix $A^{[d]}$ is used to define a $\mathfrak{g}\mathfrak{i}\mathfrak{m}$ algebra $\gimAd$. Since $\gimAd$ may possess roots with mixed signs, we quotient out by an ideal $\mfr$ that is tailor-made to capture all such roots.  The quotient algebra is called the intersection matrix  algebra and is denoted by $\imAd$.

We show that $\imAd$ is a $BC_r$-graded Lie algebra, which, in turn, allows us to invoke Allison, Benkart, and Gao's Recognition Theorem and relate $\imAd$ to an algebraic structure that is better understood [ABnG].  Combining their theorem with the knowledge that $\imAd$ is centrally closed, we conclude that, up to isomorphism, $\imAd$ is the universal covering algebra of $\sotrpoaecc$.  The algebra $\sotrpoaecc$ is like the usual matrix model $\sotrpoC$ of a finite-dimensional Lie algebra of type $B_r$, except that we now replace the field $\bbC$ with an associative algebra $\mfa$, which possesses an involution (i.e., period two antiautomorphism) $\eta$, and we involve a right $\mfa$-module $C$ that has a hermitian form $\chi: C \times C \to \mfa$.  The defining relations of the generalized intersection matrix algebra and, hence, the intersection matrix algebra, in concert with the existence of a central, graded, surjective Lie algebra homomorphism $\psi$ from $\imAd$ to $\sotrpoaecc$ allow us to understand each of $\mfa$, $\eta$, $C$, and $\chi$.  For example, we get (i) two generators of $\mfa$, namely $x$ and $\xinv$, for every long root of the form $\pm\left( \epsilon_i + \epsilon_{i+1} \right)$; and (ii) four generators of $\mfa$, namely $y$, $\yinv$, $z$, and $\zinv$, for every other type of long root that we adjoin.  We are also able to study the relations among the generators, determine the action of the involution $\eta$, and discover that $C=0$ and $\chi = 0$.  Through constructing a surjective Lie algebra homomorphism $\varphi: \gimAd \to \sotrpoaecc$ we verify that we indeed have a complete description of the ``coordinate algebra'' $\mfa$.

	Our work continues the line of research initiated by Berman and Moody, and Benkart and Zelmanov. Berman and Moody were the first to find realizations of intersection matrix algebras over Lie algebras graded by root systems of types $A_r$ ($r \geq 2$), $D_r$, $E_6$, $E_7$, and $E_8$ [BrM].  Benkart and Zelmanov found realizations of intersection matrix algebras over Lie algebras graded by root systems of types $A_1$, $B_r$, $C_r$, $F_4$, and $G_2$ [BnZ].  In this paper, we find realizations of intersection matrix algebras over Lie algebras graded by root systems of type $BC_r$ with grading subalgebras of type $B_r$ ($r \geq 3$).

\section{Multiply affinizing Cartan matrices}

In this paper, we focus on generalized intersection matrix algebras that arise from multiply affinizing a Cartan matrix of type $B_r$, where $r \geq 3$, with long roots in the root system $\Delta_{B_r}$.

	Consider a root system of type $B_r$.  Up to isomorphism, $\Delta_{B_r}$ may be described as
		\begin{displaymath}
			\Delta_{B_r} = \left\{ \pm \epsilon_i \pm \epsilon_j: \ 1 \leq i \neq j \leq r \right\} \cup \left\{ \pm \epsilon_i: \ i = 1, \ldots, r \right\}.				
		\end{displaymath}
	Once we fix an ordering of the simple roots $\alpha_1, \ldots, \alpha_r$ in a base $\Pi$, the Cartan matrix $A$ is described by
		\begin{displaymath}
			A_{ij} = \frac{2 \left( \alpha_i, \alpha_j \right)_{\text{Killing}} }{ \left( \alpha_i, \alpha_i \right)_{\text{Killing}}}, \qquad \text{for } 1 \leq i,j \leq r.
		\end{displaymath} 
	Choose any $d$ long roots in $\Delta_{B_r}$, say $\alpha_{r+1}, \ldots, \alpha_{r+d}$, and consider the $r+d$ by $r+d$ matrix $A^{[d]}$ given by
		\begin{displaymath}
			A^{[d]}_{ij} = \frac{2 \left( \alpha_i, \alpha_j \right)_{\text{Killing}} }{ \left( \alpha_i, \alpha_i \right)_{\text{Killing}}}, \qquad \text{for } 1 \leq i,j \leq r+d,
		\end{displaymath}
	with respect to the ordering $(\alpha_1, \ldots, \alpha_r, \alpha_{r+1}, \ldots, \alpha_{r+d})$ of the $r$ roots in the base $\Pi$ plus the $d$ ``adjoined'' roots.  The axioms of a root system tell us that all the entries of $A^{[d]}$ are integers.  Moreover, since the Killing form is symmetric, we have $A^{[d]}_{ji}=0$ if $A^{[d]}_{ij}=0$, or if $A^{[d]}_{ij}$ and $A^{[d]}_{ji}$ are nonzero, then they share the same sign.  In other words, $A^{[d]}$ is a generalized intersection matrix.  

 Since the ``d-affinized'' Cartan matrix $A^{[d]}$ is a generalized intersection matrix, $\gimAd$ is a generalized intersection matrix algebra.

	Note that if we affinize the Cartan matrix A of type $B_r$ with the negative of the highest long root of $\Delta_{B_r}$
then the resulting generalized intersection matrix algebra $\mathfrak{g}\mathfrak{i}\mathfrak{m}\left(A^{[1]}\right)$ is the affine Kac-Moody Lie algebra of type $B_r^{(1)}$.

\section{Intersection matrix algebras}\label{sect:im algebras}

	Fix a Cartan matrix $A$ of type $B_r$ ($r \geq 3$) with, say, $\alpha_1$, $\alpha_2$, $\ldots$, $\alpha_r$ being the simple roots in a base of $\Delta_{B_r}$ that were used to form $A$.  Let 
	\begin{itemize}
		\item
			$\Omega =$ set of all long roots of the form $\pm \big( \eipej{i}{i+1} \big) $ that we adjoin,
		\item
			$\Theta =$ set of all remaining long roots that are adjoined, 
		\item
			$N_{\mu} = $ the number of copies of the long root $\mu$ we have adjoined, and 
		\item
			$d = \sum_{\mu \in \Omega \cup \Theta} N_{\mu}$.
	\end{itemize}

	Let $A^{[d]}$ be the resulting generalized intersection matrix and $\gimAd$ the corresponding generalized intersection matrix algebra.

	We begin a move towards a quotient algebra of $\gimAd$ using a slight generalization of the work done by Benkart $\&$ Zelmanov [BnZ].  Let $\Gamma$ be the integer lattice generated by the $\Delta$, where
		\begin{displaymath}
			\Delta = \left\{ \pm \epsilon_i \pm \epsilon_j : \ 1 \leq i \neq j \leq r \right\}  \cup \left\{ \pm \epsilon_i, \ \pm 2\epsilon_i : \ i = 1, \ldots, r \right\},
		\end{displaymath}
	is a root system of type $BC_r$.

	We define a $\Gamma$-grading on $\gimAd$ as follows:
		\begin{displaymath}
			\deg e_i = \alpha_i = -\deg f_i, \quad \deg h_i = 0
		\end{displaymath}
	for $i = 1, \ldots, r$, and
		\begin{displaymath}
			\deg e_{\mu, i} = \mu = -\deg f_{\mu,i}, \quad \deg h_{\mu,i} = 0
		\end{displaymath}
	for $\mu \in \Omega \cup \Theta$ and $i = 1, \ldots, N_{\mu}$.  

	Next, we define the \emph{radical} $\mfr$ of $\gimAd$ to be the ideal generated by the root spaces $\gimAd_{\gamma}$ where $\gamma \notin \Delta \cup \{0\}$.  Since the ideal $\mfr$ is homogeneous, the resulting quotient algebra
		\begin{displaymath}
			\imAd = \gimAd / \mfr
		\end{displaymath}
	is also $\Gamma$-graded. Moreover,
		\begin{displaymath}
			\imAd_{\gamma} = 0, \quad \text{ if } \gamma \notin \Delta \cup \{0\}.
		\end{displaymath}
	We call $\imAd$ the \emph{intersection matrix ($\mathfrak{i}\mathfrak{m}$) algebra} corresponding to the generalized intersection matrix algebra $\gimAd$.






\subsection{$\imAd$ is $BC_r$-graded}\label{subsect:imAd is BCrgraded}

Allison, Benkart, and Gao gave the following definition of a Lie algebra graded by a root system of type $BC$.

\begin{defn}[{[ABnG]}]\emph{
	Let $r$ be a positive integer greater than or equal to $3$.  A Lie algebra $L$ over $\bbC$ is \emph{graded by the root system $BC_r$} or is \emph{$BC_r$-graded} with a \emph{grading subalgebra of type $B_r$} if
	\begin{enumerate}
		\item[(i)]
			$L$ contains, as a subalgebra, a finite-dimensional simple Lie algebra $\mfg$ whose root system relative to a Cartan subalgebra $\mfh = \mfg_0$ is $\Delta_{B_r}$;
		\item[(ii)]
			$L = \bigoplus_{\mu \in \Delta \cup \{0\}} L_{\mu}$, where $L_{\mu} = \left\{ x \in L | \ [h,x] = \mu(h)x \text{ for all } h \in \mfh \right\}$ for $\mu \in \Delta \cup \{0\}$, and $\Delta$ is the root system of type $BC_r$; and
		\item[(iii)]
			$L_0 = \sum_{\mu \in \Delta} \left[ L_{\mu}, L_{-\mu} \right]$.
	\end{enumerate}
}\end{defn}

	\begin{propn}\label{propn:imAd is BCrgraded}
		The algebra $\imAd$ is $BC_r$-graded with a grading subalgebra of type $B_r$.
	\end{propn}
	
		\begin{proof}
			The subalgebra in $\imAd$ generated by $e_1 + \mfr, \ldots, h_r + \mfr$, due to the relations on these elements induced by the relations on their preimages in $\gimAd$, is isomorphic to a finite-dimensional simple Lie algebra $\mfg$ of type $B_r$.  We have already shown in $\S 3.1$ that $\imAd$ is $\Gamma$-graded with $\imAd_{\gamma} = 0$ if $\gamma \notin \Delta \cup \{0\}$.  That is,
				\begin{displaymath}
					\imAd = \bigoplus_{\mu \in \Delta \cup \{0\}} \imAd_{\mu}.
				\end{displaymath}
Finally, our initial degree assignments for the generators of $\gimAd$, the $\mathfrak{g}\mathfrak{i}\mathfrak{m}$ algebra relations like $h_i = \left[ e_i, \ f_i \right]$ and $h_{\mu} = \left[ e_{\mu}, \ f_{\mu} \right]$, and the fact that movement into the $0$ root space can only occur  by bracketing an element from an $\imAd_{\mu}$ space with one from the $\imAd_{-\mu}$ space all combine to lead us to the conclusion that 
				\begin{displaymath}
					\imAd_0 = \sum_{\mu \in \Delta} \left[ \imAd_{\mu}, \ \imAd_{-\mu} \right].
				\end{displaymath}
		\end{proof}

\subsection{$\imAd$ is centrally closed}

	\begin{propn}
		The algebra $\gimAd$ is a perfect Lie algebra.
	\end{propn}
	
	\begin{proof}
			Being a Lie algebra, $\gimAd$ is closed under the operation of taking brackets; hence $ \left[ \gimAd, \ \gimAd \right] \subset \gimAd$.
			To show the reverse inclusion, it suffices to show that all of the generators of $\gimAd$ lie in $\left[ \gimAd, \ \gimAd \right]$. But this is indeed the case because the generators $e_i$, $f_i$, $h_i$ (for $1 \leq i \leq r$) and the $e_{\mu,i}$, $f_{\mu,i}$, $h_{\mu,i}$, which arise from adjoining the $i^{\text{th}}$ copy of a long root $\mu$, satisfy the relations $(R1)$ of Definition \ref{defn:gim algebra}.
			
			
		\end{proof}
	
	
	Our next theorem is Proposition 1.6 in [BnZ] adapted to our context.
	
	\begin{thm}
		The algebra $\imAd$ is centrally closed.
	\end{thm}
	\begin{proof}
		Let $\left( \tU, \phi \right)$ be the universal covering algebra of $\imAd$.  Let $\mfg$ be the simple finite dimensional subalgebra of type $B$ contained in $\imAd$ with Cartan subalgebra $\mfh$ whose root space decomposition induces a $BC$-gradation on $\imAd$.  The preimage $\phi^{-1}(\mfh)$ of $\mfh$ contains $\ker \phi$.  Since $\phi$ is a central map, $\ker \phi$ lies in the centre of $\tU$.  So 
			\begin{displaymath}
				\mfh' = \phi^{-1}(\mfh) / \ker\phi
			\end{displaymath}
		acts on $\tU$ via the adjoint action.  If $h' \in \mfh'$, $\phi(h') = h \in \mfh$, and $\mu(t) \in \bbC[t]$ is the minimal polynomial of $\text{ad}_{\UuL}(h)$, then
			\begin{displaymath}
				\mu\left( \text{ad}_{\tU}(h') \right) \left( \tU \right) \subset \ker \phi.
			\end{displaymath}
		So $\text{ad}_{\tU}(h')$ satisfies the polynomial $t\mu(t)$.  Therefore $\tU$ is a sum of root spaces with respect to $\text{ad}_{\tU} \mfh'$, and $\tU_{\gamma} \neq (0)$ if and only if $\gamma \in \Delta \cup \{0\}$.  So $\phi$ induces an isomorphism between the nonzero root spaces of $\tU$ and those of $\imAd$.  Moreover, 
			\begin{displaymath}
				\tU_0 = \sum_{\gamma \in \Delta} \left[ \tU_{-\gamma}, \tU_{\gamma} \right] + \ker \phi,
			\end{displaymath}
		implies that
			\begin{displaymath}
				\left[ \tU_0, \tU_0 \right] \subset \sum_{\gamma \in \Delta} \left[ \tU_{-\gamma}, \tU_{\gamma} \right].
			\end{displaymath}
		Since $\tU = \left[ \tU, \tU \right]$, it follows that
			\begin{displaymath}
				\tU_0 = \left[ \tU_0, \tU_0 \right] + \sum_{\gamma \in \Delta} \left[ \tU_{-\gamma}, \tU_{\gamma} \right] = \sum_{\gamma \in \Delta} \left[ \tU_{-\gamma}, \tU_{\gamma} \right].
			\end{displaymath}
		Consequently, $\phi$ is an isomorphism.
\end{proof}

\section{Recognition Theorem}\label{sect:recognition theorem}

	The following construction, given in Example 1.23 of [ABnG], is a more general version of the classical construction of $\sotrpoC$, the simple Lie algebra of type $B_r$.

	Let $r$ be a positive integer; $\mfa$ be a unital associative algebra over $\bbC$ with an involution (i.e., period two anti-automorphism) $\eta$; $C$ be a right $\mfa$-module with a hermitian form $\chi: C \times C \to \mfa$, i.e. a biadditive map $\chi: C \times C \to \mfa$ satisfying
		\begin{itemize}
			\item[]
				$\chi(c,c'\cdot a) = \chi(c,c')\cdot a$
			\item[]
				$\chi(c\cdot a, c') = \eta(a) \cdot \chi(c,c')$
			\item[]
				$\chi(c,c') = \eta \left( \chi(c',c) \right)$
		\end{itemize}
	for $c,c' \in C$, $a \in \mfa$; and	$G$ be the $(2r+1) \times (2r+1)$ matrix
	\begin{displaymath}
		G 
			= \left[
				\begin{array}{ccccc}
					0 & 0 & \cdots & 0 & 1 \\
					0 & 0 & \cdots & 1 & 0 \\
					\vdots & \vdots & \ddots & \vdots & \vdots \\
					0 & 1 & \cdots & 0 & 0 \\
					1 & 0 & \cdots & 0 & 0 \\
				\end{array}
			\right].
	\end{displaymath}
	
	Also, given any $c \in C$, define $\chi_c \in C^*$ by
		\begin{displaymath}
			\chi_c(c') := \chi(c,c'),
		\end{displaymath}	
	for any $c' \in C$, and	given any $\usc = \left[ \begin{array}{c} c_1 \\ \vdots \\ c_n \end{array} \right] \in C^{2r+1}$, define
		\begin{displaymath}
			\chi_{\usc} := \left[ \begin{array}{c} \chi_{c_1} \\ \vdots \\ \chi_{c_n} \end{array} \right] \in \left( C^* \right)^{2r+1} .
		\end{displaymath}		

	Now let
		\begin{align*}
			\mfA(\chi) := \Bigg\{ N \in \EndaC : \chi(Nc,c') & + \chi(c,Nc') = 0 \text{ for all } c,c' \in C   \Bigg\},
		\end{align*} 
	and
		\begin{displaymath}
			\mfA := \left\{  \left[ \begin{array}{cc} M & \chi_{\usc} \\ \usc^t G & N  \end{array} \right] : M \in \text{M}_{2r+1}(\mathfrak{a}), \ (M^{\eta})^t G + GM = 0, \ \usc \in C^{2r+1}, \ N \in \mfA(\chi)   \right\}.\\ 
		\end{displaymath} 

	It can be checked that $\mfA$ is a Lie algebra that contains a simple Lie algebra 
		\begin{displaymath}
			\mfg = \left\{ \left[ \begin{array}{cc} M & 0 \\ 0 & 0 \end{array} \right] : M \in \text{M}_{2r+1}(\bbC), \ M^t G + GM = 0   \right\},
		\end{displaymath}
	of type $B_r$.  If $\mfh$ denotes the Cartan subalgebra of diagonal matrices in $\mfg$, then the adjoint action of $\mfh$ on $\mfA$ induces a root space decomposition
		\begin{displaymath}
			\mfA = \bigoplus_{\mu \in \Delta \cup \{0\}} \mfA_{\mu},
		\end{displaymath}
	where
		\begin{displaymath}
			\mfA_{\mu} = \left\{ T \in \mfA : [h,T] = \mu(h) \, T \text{  for all } h \in \mfh \right\}.
		\end{displaymath}   
	The following abbreviated notation helps describe these root spaces:  for $v = \left[ \begin{array}{c} v_1 \\ \vdots \\ v_{2r+1} \end{array} \right] \in \bbC^{2r+1}$ and $c \in C$, let $vc := \left[ \begin{array}{c} v_1 c \\ \vdots \\ v_{2r+1} c \end{array} \right] \in C^{2r+1}$.  Then $C^{2r+1} = \bigoplus_{i=1}^{2r+1} e_i C$, where $e_1, \ldots, e_{2r+1}$ is the standard basis for $\bbC^{2r+1}$.  Letting $B$ denoting the set of skew-symmetric elements of $\mfa$ relative to the involution $\eta$, we have
		\begin{align*}
			& \mfA_{\eimej{i}{j}} = \left\{ E_{i,j}(a) + E_{2r+2-j,2r+2-i}(-\eta(a)): \ a \in \mfa \right\}, \quad 1 \leq i \neq j \leq r, \\
			& \mfA_{\eipej{i}{j}} = \left\{ E_{i,2r+2-j}(a) + E_{j,2r+2-i}(-\eta(a)): \ a \in \mfa \right\}, \quad 1 \leq i,j \leq r, \\
			& \mfA_{\meimej{i}{j}} = \left\{ E_{2r+2-i,j}(a) + E_{2r+2-j,i}(-\eta(a)): \ a \in \mfa \right\}, \quad 1 \leq i,j \leq r, \\				& \mfA_{\epsilon_i} = \left\{ \left[ \begin{array}{cc} 0 & \chi_{e_i c} \\ (e_{2r+2-i} c)^t & 0 \end{array} \right]: \ c \in C \right\} \\
			& \qquad + \left\{ E_{i,r+1}(a) + E_{r+1,2r+2-i}(-\eta(a)): \ a \in \mfa \right\}, \quad 1 \leq i \leq r, \\					
			& \mfA_{-\epsilon_i} = \left\{ \left[ \begin{array}{cc} 0 & \chi_{e_{2r+2-i} c} \\ (e_i c)^t & 0 \end{array} \right]: \ c \in C \right\} \\
			& \qquad + \left\{ E_{r+1,i}(a) + E_{2r+2-i,r+1}(-\eta(a)): \ a \in \mfa \right\}, \quad 1 \leq i \leq r, \\								
			& \mfA_0 = \left\{ \sum_{i=1}^r E_{ii}(a) + E_{2r+2-i,2r+2-i}(-\eta(a)): \ a \in \mfa \right\} \\			
			& \qquad + \left\{ \left[ \begin{array}{cc} 0 & 0 \\ 0 & N \end{array} \right] : \ N \in \mfA(\chi) \right\} + \left\{ E_{r+1,r+1}(b): \ b \in B \right\} \\
			& \qquad \qquad \left\{ \left[ \begin{array}{cc} 0 & \chi_{e_{r+1} c} \\ (e_{r+1} c)^t & 0 \end{array} \right] : \ c \in C \right\}.
		\end{align*}

	The subalgebra 
		\begin{displaymath}
			\sotnpoaecc := \sum_{\mu \in \Delta} \mfA_{\mu} \ + \ \sum_{\mu \in \Delta} \left[ \mfA_{\mu}  , \,  \mfA_{-\mu}  \right]
			\end{displaymath}
	of $\mfA$ has the root spaces
		\begin{displaymath}
			\sotrpoaecc_0 = \sotrpoaecc \cap \mfA_0
		\end{displaymath}
	and 
		\begin{displaymath}
			\sotrpoaecc_{\mu} = \mfA_{\mu} \text{ for } \mu \in \Delta.
		\end{displaymath} 
	In particular,
		\begin{displaymath}
			\sotnpoaecc_0 = \sum_{\mu \in \Delta} \left[ \sotnpoaecc_{\mu}, \, \sotnpoaecc_{-\mu} \right].
		\end{displaymath}

\noindent \emph{\bf Remark:}  [ABnG] use the notation $L$ to refer to the Lie algebra that we are calling $\sotnpoaecc$.

	
	 We use the following notation to shorten the description of elements in $\sotnpoaecc$:  Given any $1 \leq k \leq r$ and $a \in \mfa$, let
	\begin{align*}
		& \Evert{k}(a) := E_{k,r+1}(a) + E_{r+1,2r+2-k}(-\eta(a)), \\
		&	\Ehort{k}(a) := E_{r+1,k}(a) + E_{2r+2-k,r+1}(-\eta(a)),
	\end{align*}
and for any $1 \leq p,q \leq r$ and $a \in \mfa$, let
	\begin{align*}
		& \EUL{p}{q}(a) := E_{p,q}(a) + E_{2r+2-q,2r+2-p}(-\eta(a)), \\
		&	\EUR{p}{q}(a) := E_{p,2r+2-q}(a) + E_{q,2r+2-p}(-\eta(a)), \\
		&	\EBL{p}{q}(a) := E_{2r+2-p,q}(a) + E_{2r+2-q,p}(-\eta(a)).
	\end{align*}

We often also denote the involution $\eta$ on $\mfa$ by $\bar{\cdot}$.  So, for example, we would write $\EBL{p}{q}(a)$ above as $E_{2r+2-p,q}(a) + E_{2r+2-q,p}(-\ovla).$

Allison, Benkart, and Gao's classification results on $BC_r$-graded Lie algebras say the following in our setting [ABnG]:

	\begin{thm}[{[ABnG], Thm. 3.10}]
	Let $r \geq 3$.  Then $L$ is $BC_r$-graded with grading subalgebra $\mfg$ of type $B_r$ if and only if there exists an associative algebra $\mfa$ with involution $\eta$, and an $\mfa$-module $C$ with a hermitian form $\chi$	such that $L$ is centrally isogenous to the $BC_r$-graded Lie algebra $\sotnpoaecc$.	
		\end{thm}

Since $\imAd$ is $BC_r$-graded with a grading subalgebra of type $B_r$ and is centrally closed, we have the following result.


		\begin{cor}\label{cor:imAd is uca of sotrpo}
		 The intersection matrix algebra $\imAd$ is isomorphic to the universal covering algebra of the Lie algebra $\sotrpoaecc$.  In particular, there exists a graded, central, surjective Lie algebra homomorphism $\psi: \imAd \to \sotnpoaecc$.
		\end{cor}










\section{Arriving at a ``minimal'' understanding of $\mfa$, $\eta$, $C$, and $\chi$}\label{sect:minimal understanding}

The graded nature of the map $\psi: \imAd \to \sotrpoaecc$ along with the relations among the generating elements of $\imAd$ allow us to study each of components $\mfa$, $\eta$, $C$, and $\chi$ involved in $\sotrpoaecc$.

		Since the elements $e_1 + \mfr, \ldots, e_r + \mfr$, $f_1 + \mfr, \ldots, f_r + \mfr$, $h_1 + \mfr, \ldots, h_r + \mfr$ in $\imAd$ generate a subalgebra isomorphic to a simple Lie algebra of type $B_r$, and since $\psi$ is a graded homomorphism, we may assume without loss of generality that (after relabeling the $e_i + \mfr$, $f_i + \mfr$, and $h_i + \mfr$ as $e_i$, $f_i$, and $h_i$, respectively)
			\begin{align*}
				& \psi(e_i) = \EUL{i}{i+1}(1), \quad \text{ for } 1 \leq i \leq r-1, \\
				& \psi(e_r) = \Evert{r}\left( \sqrt{2} \right), \\
				& \psi(f_i) = \EUL{i+1}{i}(1), \quad \text{ for } 1 \leq i \leq r-1, \\
				& \psi(f_r) = \Ehort{r}\left( \sqrt{2} \right), \\
				& \psi(h_i) = \EUL{i}{i}(1) + \EUL{i+1}{i+1}(-1), \ \text{ for } 1 \leq i \leq r-1, \text{ and} \\
				& \psi(h_r) = \EUL{r}{r}(2).
			\end{align*}

\noindent \emph{\bf Remark:}  Here we are using the notation established in Section \ref{sect:recognition theorem}.  The generators of $\imAd$ coming from a simple root $\alpha_j \in \Pi$ are denoted by $e_j$, $f_j$, and $h_j$, while the generators coming from an $i^{\text{th}}$ copy of an adjoined root $\alpha \in \Delta_{B_r}$ are denoted by $e_{\alpha,i}$, $f_{\alpha,i}$, and $h_{\alpha,si}$.


\subsection{Understanding the invertibility of some coordinates of $\mfa$}




\begin{propn}
	Let $e_{\eimej{p}{q},i}$, $f_{\eimej{p}{q},i}$, $h_{\eimej{p}{q},i}$ be the generators of $\imAd$ that result from adjoining the $i^{\text{th}}$ copy of a long root $\eimej{p}{q}$ ($1 \leq p,q \leq r$, $p \neq q$).  If $\psi\left( e_{\eimej{p}{q},i} \right) = \EUL{p}{q}( a )$ for some $a \in \mfa$, then $a$ is an invertible element and $\psi\left( f_{\eimej{p}{q},i} \right) = \EUL{q}{p}( a^{-1} )$.
\end{propn}

\begin{proof}
		Since $\psi$ is a graded homomorphism, 
			\begin{align*}
				& \psi\left( e_{\eimej{p}{q},i} \right) = \EUL{p}{q}( a ), \text{ and} \\
				& \psi\left( f_{\eimej{p}{q},i} \right) = \EUL{q}{p}( a' ),
			\end{align*}
		for some $a$, $a'$ in the associative algebra $\mfa$.  
		
		Assume that $p < q$. (If $q < p$, then we reverse the roles of $p$ and $q$ in the following analysis.)  Then
			\begin{align*}
				& \left[ \left[ \psi\left( e_{\eimej{p}{q},i} \right), \ \psi\left( f_{\eimej{p}{q},i} \right) \right], \ \psi\left( e_q \right) \right] \\
				& 
				= \left\{
						\begin{array}{ll}
							\left[  \left[ \EUL{p}{q}(a), \ \EUL{q}{p}(a') \right], \ \EUL{q}{q+1}(1) \right] & \text{ if } q < r \\
							\left[  \left[ \EUL{p}{q}(a), \ \EUL{q}{p}( a' ) \right], \ \Evert{r}\left( \sqrt{2} \right) \right] & \text{ if } q = r 
						\end{array}
					\right. \\
				& 
				= \left\{
						\begin{array}{ll}
							\left[ \EUL{p}{p}( a a' ) + \EUL{q}{q}( -a' a ), \ \EUL{q}{q+1}(1)  \right] & \text{ if } q < r \\
							\left[ \EUL{p}{p}( a a' ) + \EUL{q}{q}( -a' a ), \ \Evert{r}\left( \sqrt{2} \right) \right] & \text{ if } q = r 
						\end{array}
					\right. \\
				&
				= \left\{ 
						\begin{array}{ll}
							\EUL{q}{q+1}( -a' a ) & \text{ if } q < r \\
							\Evert{r}\left( -\sqrt{2} \ a' a \right) & \text{ if } q = r
						\end{array}
					\right.
			\end{align*}
		But 
			\begin{displaymath}
				\left[ e_{\eimej{p}{q},i}, \ f_{\eimej{p}{q},i} \right] = h_{\eimej{p}{q},i}
			\end{displaymath}
		and 
			\begin{align*}
				\left[ h_{\eimej{p}{q},i}, \ e_q \right] 
				& = \left\{
							\begin{array}{ll}
								A_{\eimej{p}{q}, \eimej{q}{q+1}} \ e_q & \text{ if } q < r \\
								A_{\eimej{p}{q}, \epsilon_r} \ e_q & \text{ if } q = r 
							\end{array}
						\right. \\
				& = - e_q
			\end{align*}
		imply that
			\begin{align*}
				& \left[ \left[ \psi\left( e_{\eimej{p}{q},i} \right), \ \psi\left( f_{\eimej{p}{q},i} \right) \right], \ \psi\left( e_q \right) \right] \\
				& = \left[ \psi\left( h_{\eimej{p}{q},i} \right), \ \psi\left( e_q \right) \right] \\
				& = -\psi\left( e_q \right) \\
				& = \left\{
							\begin{array}{ll}
								\EUL{q}{q+1}(-1) & \text{ if } q < r \\
								\Evert{r}\left( -\sqrt{2} \right)& \text{ if } q=r
							\end{array}
						\right.
			\end{align*}
		Hence, in either case, we get that 
			\begin{equation}\label{eq:eimej inv}
				a' a = 1.
			\end{equation}
			
		Continuing forward, on the one hand, we have
			\begin{align*}
				& \left[ \left[ \psi\left( e_{\eimej{p}{q},i} \right), \ \psi\left( f_{\eimej{p}{q},i} \right) \right], \ \psi\left( e_p \right) \right] \\
				& =	\left[  \left[ \EUL{p}{q}( a ), \ \EUL{q}{p}( a' ) \right], \ \EUL{p}{p+1}(1) \right] \\
				& = \left[ \EUL{p}{p}( a a' ) + \EUL{q}{q}( -a' a ), \ \EUL{p}{p+1}(1)  \right] \\
				&
				= \left\{ 
						\begin{array}{ll}
							\EUL{p}{p+1}( a a' ) & \text{ if } q-p \geq 2 \\
							\EUL{p}{p+1}( a a' + a' a ) & \text{ if } q = p+1
						\end{array}
					\right.
			\end{align*}
		But since
			\begin{align*}
				& \left[ \left[ e_{\eimej{p}{q},i}, \ f_{\eimej{p}{q},i} \right], \ e_p \right] \\
				& = \left[ h_{\eimej{p}{q},i}, \ e_p \right] \\
				& = A_{\eimej{p}{q}, \eimej{p}{p+1}} \ e_p \\
				& = \left( 1 + \delta_{q,p+1} \right) \ e_p\\
				& = 
					\left\{
						\begin{array}{ll}
							e_p & \text{ if } q \geq p+2 \\
							2 e_p & \text{ if } q = p+1
						\end{array}
					\right.
			\end{align*}
		we have
			\begin{displaymath}
				\left[ \left[ \psi\left( e_{\eimej{p}{q},i} \right), \ \psi\left( f_{\eimej{p}{q},i} \right) \right], \ \psi\left( e_p \right) \right]
				= \left\{ 
						\begin{array}{ll}
							\EUL{p}{p+1}(1) & \text{ if } q \geq p+2 \\
							\EUL{p}{p+1}(2) & \text{ if } q = p+1
						\end{array}
					\right.
			\end{displaymath}		
		So if $q \geq p+2$, then $ a a' = 1$, and if $q=p+1$, then $ a a' + a' a = 2$.  But we already know by equation (\ref{eq:eimej inv}) that $a' a = 1$.  Hence, in either case,
			\begin{displaymath}
				a a' = 1.
			\end{displaymath}
		So
			\begin{displaymath}
				 a' = a^{-1}, 
			\end{displaymath}
		and if $\psi\left( e_{\eimej{p}{q},i} \right) = \EUL{p}{q}( a )$, then $\psi\left( f_{\eimej{p}{q},i} \right) = \EUL{q}{p}( a^{-1} )$.
\end{proof}

		


\begin{propn}
	Let $e_{\eipej{p}{q},i}$, $f_{\eipej{p}{q},i}$, $h_{\eipej{p}{q},i}$ be the generators of $\imAd$ that result from adjoining the $i^{\text{th}}$ copy of a long root $\eipej{p}{q}$ ($1 \leq p,q \leq r$, $p \neq q$).  If $\psi\left( e_{\eipej{p}{q},i} \right) = \EUR{p}{q}( b )$ for some $b \in \mfa$, then $b$ is an invertible element and $\psi\left( f_{\eipej{p}{q},i} \right) = \EBL{q}{p}( b^{-1} )$.
\end{propn}

\begin{proof}
Since $\psi$ is a graded homomorphism, 
	\begin{align*}
		& \psi\left( e_{\eipej{p}{q},i} \right) = \EUR{p}{q}( b ), \text{ and} \\
		& \psi\left( f_{\eipej{p}{q},i} \right) = \EBL{q}{p}( b' ),
	\end{align*}
for some $b$, $b'$ in the associative algebra $\mfa$.  
		
Without loss of generality, we may assume that $p < q$.  Then
	\begin{align*}
		& \left[ \left[ \psi\left( e_{\eipej{p}{q},i} \right), \ \psi\left( f_{\eipej{p}{q},i} \right) \right], \ \psi\left( e_q \right) \right] \\
		& 
		= \left\{
				\begin{array}{ll}
					\left[ \left[ \EUR{p}{q}( b ), \ \EBL{q}{p}( b' ) \right], \ \EUL{q}{q+1}(1) \right] & \text{ if } q < r \\
					\left[ \left[ \EUR{p}{q}( b ), \ \EBL{q}{p}( b' ) \right], \ \Evert{r}\left( \sqrt{2} \right) \right] & \text{ if } q = r 
				\end{array}
			\right. \\	
		& 
		= \left\{
				\begin{array}{ll}
					\left[ \EUL{p}{p}( b b' ) + \EUL{q}{q}\left( \eta(b ) \, \eta(b' ) \right), \ \EUL{q}{q+1}(1) \right] & \text{ if } q < r \\
					\left[ \EUL{p}{p}( b b' ) + \EUL{q}{q}\left( \eta( b  ) \, \eta( b' ) \right), \ \Evert{r}\left( \sqrt{2} \right) \right] & \text{ if } q = r 
				\end{array}
			\right. \\						
		&
		= \left\{ 
				\begin{array}{ll}
					\EUL{q}{q+1}\left( \eta( b ) \, \eta( b' ) \right) & \text{ if } q < r \\
					\Evert{r}\left( \sqrt{2} \ \eta( b ) \, \eta( b' ) \right) & \text{ if } q = r
				\end{array}
			\right.
	\end{align*}
But 
	\begin{displaymath}
		\left[ e_{\eipej{p}{q},i}, \ f_{\eipej{p}{q},i} \right] = h_{\eipej{p}{q},i}
	\end{displaymath}
and 
	\begin{align*}
		\left[ h_{\eipej{p}{q},i}, \ e_q \right] 
		&	= \left\{
					\begin{array}{ll}
						A_{\eipej{p}{q}, \eimej{q}{q+1}} \ e_q & \text{ if } q < r \\
						A_{\eipej{p}{q}, \epsilon_r} \ e_q & \text{ if } q = r 
					\end{array}
				\right. \\
		& = e_q
	\end{align*}
imply that
	\begin{align*}
		& \left[ \left[ \psi\left( e_{\eipej{p}{q},i} \right), \ \psi\left( f_{\eipej{p}{q},i} \right) \right], \ \psi\left( e_q \right) \right] \\
		& = \left[ \psi\left( h_{\eipej{p}{q},i} \right), \ \psi\left( e_q \right) \right] \\
		& = \psi\left( e_q \right) \\
		& = \left\{
					\begin{array}{ll}
						\EUL{q}{q+1}(1) & \text{ if } q < r \\
						\Evert{r}\left( \sqrt{2} \right)& \text{ if } q=r
					\end{array}
				\right.
	\end{align*}
whence
	\begin{displaymath}
		\eta( b ) \, \eta( b' ) = 1.
	\end{displaymath}
Taking $\eta$ of both sides and using the fact that it is an anti-automorphism of order $2$, we get that
	\begin{equation}\label{eq:eipej inv}
		b' b = 1.
	\end{equation}

We also have
	\begin{align*}
		& \left[ \left[ \psi\left( e_{\eipej{p}{q},i} \right), \ \psi\left( f_{\eipej{p}{q},i} \right) \right], \ \psi\left( e_p \right) \right] \\
		& 
			= \left[ \left[ \EUR{p}{q}( b ), \ \EBL{q}{p}( b' ) \right], \ \EUL{p}{p+1}(1) \right] \\
		& 
			= \left[ \EUL{p}{p}( b b' ) + \EUL{q}{q}\left( \eta( b ) \, \eta( b' ) \right), \ \EUL{p}{p+1}(1) \right] \\
		&
			= \left\{ 
				\begin{array}{ll}
					\EUL{p}{p+1}( b b' ) & \text{ if } q \geq p+2 \\
					\EUL{p}{p+1}( b b' - \eta( b ) \, \eta( b' ) ) & \text{ if } q = p+1
				\end{array}
			\right.
	\end{align*}

But since
	\begin{align*}
		& \left[ \left[ e_{\eipej{p}{q},i}, \ f_{\eipej{p}{q},i} \right], \ e_p \right] \\
		& 
			= \left[ h_{\eipej{p}{q},i}, \ e_p \right] \\
		&
			= A_{\eipej{p}{q}, \eimej{p}{p+1}} \ e_p \\
		&
			= \left( 1 - \delta_{q,p+1} \right) \ e_p \\
		&
			= \left\{
					\begin{array}{ll}
						e_p & \text{ if } q \geq p+2 \\
						0 & \text{ if } q = p+1 
					\end{array}
				\right.	
	\end{align*}

we have
	\begin{displaymath}
		\left[ \left[ \psi\left( e_{\eipej{p}{q},i} \right), \ \psi\left( f_{\eipej{p}{q},i} \right) \right], \ \psi\left( e_p \right) \right]
		= \left\{ 
				\begin{array}{ll}
					\EUL{p}{p+1}(1) & \text{ if } q-p \geq 2 \\
					0 & \text{ if } q = p+1
				\end{array}
			\right.
	\end{displaymath}		
So if $q \geq p+2$, then $b b' = 1$, and if $q=p+1$, then 
	\begin{align*}
		& b b' - \eta(b) \, \eta(b') = 0 \\
		& \implies b b' = \eta(b) \, \eta(b') \\
		& \implies b b' = \eta(b' b) \\
		& \implies b b' = \eta\left( 1 \right), \ \text{ since } b' b = 1 \text{ by eqn. } (\ref{eq:eipej inv}) \\
		& \implies b b' = 1 
	\end{align*}
In either case,
	\begin{displaymath}
		b b' = 1.
	\end{displaymath}

So, in light of the equalities $b' b = 1$ and $ b b' = 1$, we conclude that
	\begin{displaymath}
		b' = b^{-1}.
	\end{displaymath}
In particular, if $\psi\left( e_{\eipej{p}{q},i} \right) = \EUR{p}{q}( b )$, then $\psi\left( f_{\eipej{p}{q},i} \right) = \EBL{q}{p}( b^{-1} )$.						
\end{proof}




Using similar calculations as above, the following also holds.

\begin{propn}
	Let $e_{\meimej{p}{q},i}$, $f_{\meimej{p}{q},i}$, $h_{\meimej{p}{q},i}$ be the generators of $\imAd$ that result from adjoining the $i^{\text{th}}$ copy of a long root $\meimej{p}{q}$ ($1 \leq p,q \leq r$, $p \neq q$).  If $\psi\left( e_{\meimej{p}{q},i} \right) = \EBL{p}{q}( c )$ for some $c \in \mfa$, then $c$ is an invertible element and $\psi\left( f_{\meimej{p}{q},i} \right) = \EUR{q}{p}( c^{-1}  )$.
\end{propn}

\subsection{Understanding the involution $\eta$ on $\mfa$}


\begin{propn}
If $\psi\left( e_{\eipej{p}{p+1},i} \right) = \EUR{p}{(p+1)}(a)$ for some $1 \leq p \leq r-1$ and $a \in \mfa$, then $\eta (a) = a $.
\end{propn}

\begin{proof}
Observe that 
	\begin{align*}
		& \left[ \psi\left( e_{\eipej{p}{p+1},i} \right), \ \psi\left( e_p \right) \right] \\
		& = \left[ \EUR{p}{(p+1)}(a), \ \EUL{p}{p+1}(1) \right] \\
		& = \EUR{p}{p}(a)
	\end{align*}
	
But
	\begin{displaymath}
		A_{\eipej{p}{p+1}, \eimej{p}{p+1}} = 0
	\end{displaymath}
implies that
	\begin{displaymath}
		\left( \ad e_{\eipej{p}{p+1},i} \right)^{-0+1} \ e_p = \left[ e_{\eipej{p}{p+1},i}, \ e_p \right] = 0,
	\end{displaymath}
which, in turn, implies that
	\begin{displaymath}
		\left[ \psi\left( e_{\eipej{p}{p+1},i} \right), \ \psi\left( e_p \right) \right] = 0
	\end{displaymath}

Hence
	\begin{align*}
		& \EUR{p}{p}(a)  = 0 \\
		& \implies E_{p,2r+2-p} \left( -a + \eta (a) \right) = 0 \\
		& \implies -a + \eta (a) = 0
	\end{align*}
Thus
	\begin{equation}\label{eq:eppeppo invln}
			\eta (a) = a.
	\end{equation}
\end{proof}



Similarly,
\begin{propn}
If $\psi\left( e_{\meimej{p}{p+1},i} \right) = \EBL{p}{p+1}(b)$ for some $1 \leq p \leq r-1$ and $b \in \mfa$, then $\eta (b) = b $.
\end{propn}

	
	

\subsection{Understanding the relations on generators of $\mfa$}\label{subsect:relations}



\begin{propn}
If, as a consequence of adjoining an $i^{\text{th}}$ copy of the long root $\eimej{p}{q}$ and a $j^{\text{th}}$ copy of the long root $\eipej{p}{q}$, where $1 \leq p, q \leq r$ with $p \neq q$, 
	\begin{displaymath}
		\psi\!\left( e_{\eimej{p}{q},i} \right) = \EUL{p}{q}(s) \text{ and } \		\psi\!\left( e_{\eipej{p}{q},j} \right) = \EUR{p}{q}( t ),
	\end{displaymath}
for some $s,t \in \mfa$,	then
	\begin{enumerate}
		\item[(a)]
			if $|p-q| = 1$, the elements $ s $, $ t $, and $\eta(s)$ in $\mfa$ satisfy the relation
	\begin{displaymath}
		s \cdot t = t \cdot \eta( s ) ;
	\end{displaymath}		
		\item[(b)]
			if $|p-q| \geq 2$, the elements $ s, \, t,\,  \eta( s )$, and  $\eta( t )$ in $\mfa$ satisfy the relation
	\begin{displaymath}
		s \cdot \eta( t ) = t \cdot \eta( s ).
	\end{displaymath}		
	\end{enumerate}
\end{propn}

\begin{proof}
Observe that
	\begin{align*}
		& \left[ \psi\!\left( e_{\eimej{p}{q},i} \right), \ \psi\!\left( e_{\eipej{p}{q},j} \right) \right] \\
		& = \left[ \EUL{p}{q}( s ), \ \EUR{p}{q}( t ) \right] \\
		& = \EUR{p}{p}\left( -s \cdot \eta( t ) \right) \\
		& = E_{p,2r+2-p}\left( -s \cdot \eta( t ) + t \cdot \eta( s ) \right) \\
		& = 
			\left\{
				\begin{array}{ll}
					E_{p,2r+2-p}\left( -s \cdot t + t \cdot \eta( s ) \right) & \text{ if } |p-q| = 1 \\
					E_{p,2r+2-p}\left( -s \cdot \eta( t ) + t \cdot \eta( s ) \right) & \text{ if } |p-q| \geq 2
				\end{array}
			\right.
	\end{align*}
(The division into two cases in the last step follows from the use of equation (\ref{eq:eppeppo invln}).)
But since
	\begin{displaymath}
		A_{\eimej{p}{q}, \eipej{p}{q}} = 0
	\end{displaymath}
the generalized intersection matrix algebra relations tell us that 
	\begin{displaymath}
		\left( \ad e_{\eimej{p}{q},i} \right)^{-0+1} \ e_{\eipej{p}{q},j} = 0.
	\end{displaymath}
That is, 
	\begin{displaymath}
		\left[ e_{\eimej{p}{q},i}, \ e_{\eipej{p}{q},j} \right] = 0.
	\end{displaymath}
So we must have that 
	\begin{displaymath}
		\left[ \psi\!\left( e_{\eimej{p}{q},i} \right), \ \psi\!\left( e_{\eipej{p}{q},j} \right) \right] = 0.
	\end{displaymath}
This implies that
	\begin{displaymath}
		-s \cdot t + t \cdot \eta( s ) = 0, \quad \text{ if } |p-q| = 1,
	\end{displaymath}
and
	\begin{displaymath}
		-s \cdot \eta( t ) + t \cdot \eta( s ) = 0, \quad \text{ if } |p-q| \geq 2.
	\end{displaymath}
\end{proof}

Similarly,
\begin{propn}
If, as a consequence of adjoining an $i^{\text{th}}$ copy of the long root $\eimej{p}{q}$ and a $j^{\text{th}}$ copy of the long root $\meimej{p}{q}$, where $1 \leq p, q \leq r$ with $p \neq q$, 

	\begin{displaymath}
		\psi\!\left( e_{\eimej{p}{q},i} \right) = \EUL{p}{q}(s) \text{ and } \				\psi\!\left( e_{\meimej{p}{q},j} \right) = \EBL{p}{q}( t ),
	\end{displaymath}
for some $s,t \in \mfa$,	then
	\begin{enumerate}
		\item[(a)]
			if $|p-q| = 1$, the elements $ s $, $ t $, and $\eta(s)$ in $\mfa$ satisfy the relation
	\begin{displaymath}
		\eta(s) \cdot t = t \cdot s ;
	\end{displaymath}		
		\item[(b)]
			if $|p-q| \geq 2$, the elements $ s, \, t,\,  \eta( s )$, and  $\eta( t )$ in $\mfa$ satisfy the relation
	\begin{displaymath}
		\eta(s) \cdot t = \eta(t) \cdot s.
	\end{displaymath}		
	\end{enumerate}
\end{propn}


	\subsection{A description of the module $C$}

		Since $\psi$ is a graded, surjective homomorphism from $\imAd$ to $\sotnpoaecc$ and we are only adjoining long roots, we can examine the image of $\imAd$ under $\psi$ to help us understand $C$.  In order to do so, we first examine how two typical image elements in $\sotnpoaecc$ bracket with each other.  
		
		\emph{\bf Remark.}  In what follows, given any $a \in \mfa$, the notation $\ovla$ stands for $\eta(a)$.


\begin{lemma*}[{[\vertbox, \ \vertbox]}]\label{vert vert}  
Given any $1 \leq k,p \leq r$ and $a,b \in \mfa$
	\begin{displaymath}
	[ \Evert{k}(a), \ \Evert{p}(b) ] = \EUR{k}{p}(-a\,\ovlb).
	\end{displaymath}
\end{lemma*}
\begin{proof}
	\begin{align*}
		& [ \Evert{k}(a), \ \Evert{p}(b) ] \\
		& = \left[  \EPLUSR{k}{a}{-\ovla}, \ \EPLUSR{p}{b}{-\ovlb}  \right] \\
		& =  \EURquad{k}{p}{-a \, \ovlb}{p}{k}{b \, \ovla} \\
		& = \EUR{k}{p} \left( -a \, \ovlb \right).
	\end{align*}
\end{proof}

Similarly, 

\begin{lemma*}
	\begin{enumerate}
		\item[]
		\item[]\label{vert hort} 
			$[\vertbox, \ \hortbox]$
				Given any $1 \leq k,p \leq r$ and $a,b \in \mfa$
					\begin{displaymath}
						[ \Evert{k}(a), \ \Ehort{p}(b) ] = \EUL{k}{p}(ab) + \delta_{k,p} E_{r+1,r+1}\left(-ba + \ovlab\,\right).
					\end{displaymath}

		\item[]\label{vert ul}
			$[\vertbox, \ \ulbox]$ 
				Given any $1 \leq k,p,q \leq r$ and $a,b \in \mfa$
					\begin{displaymath}
						[ \Evert{k}(a), \ \EUL{p}{q}(b) ] = \delta_{k,q} \Evert{p}(-ba).
					\end{displaymath}

		\item[]\label{vert ur}
			$[\vertbox, \ \urbox\phantom{}]$
				Given any $1 \leq k,p,q \leq r$ and $a,b \in \mfa$
					\begin{displaymath}
						[ \Evert{k}(a), \ \EUR{p}{q}(b) ] = 0.
					\end{displaymath}

		\item[]\label{vert bl}
			$[\vertbox, \ \blbox]$
				Given any $1 \leq k,p,q \leq r$ and $a,b \in \mfa$
					\begin{displaymath}
						[ \Evert{k}(a), \ \EBL{p}{q}(b) ] = \delta_{k,p} \Ehort{q}(-\ovla\,b) + \delta_{k,q} \Ehort{p}(\ovlab\,).
					\end{displaymath}

	\item[]\label{hort hort}
		$[\hortbox, \ \hortbox]$
			Given any $1 \leq k,p \leq r$ and $a,b \in \mfa$
				\begin{displaymath}
					[ \Ehort{k}(a), \ \Ehort{p}(b) ] = \EBL{k}{p}(-\ovla\,b).
				\end{displaymath}

		\item[]\label{hort ul}
			$[\hortbox, \ \ulbox]$
				Given any $1 \leq k,p,q \leq r$ and $a,b \in \mfa$
					\begin{displaymath}
						[ \Ehort{k}(a), \ \EUL{p}{q}(b) ] = \delta_{k,p} \Ehort{q}(ab).
					\end{displaymath}

		\item[]\label{hort ur}
			$[\hortbox, \ \urbox\phantom{}]$
				Given any $1 \leq k,p,q \leq r$ and $a,b \in \mfa$
					\begin{displaymath}
						[ \Ehort{k}(a), \ \EUR{p}{q}(b) ] = \delta_{k,p} \Evert{q}(-\ovlba\,) + \delta_{k,q} \Evert{p}(b\ovla).
					\end{displaymath}

		\item[]\label{hort bl}
			$[\hortbox, \ \blbox]$
				Given any $1 \leq k,p,q \leq r$ and $a,b \in \mfa$
					\begin{displaymath}
						[ \Ehort{k}(a), \ \EBL{p}{q}(b) ] = 0.
					\end{displaymath}

		\item[]\label{ul ul}
			$[\ulbox, \ \ulbox]$
				Given any $1 \leq p,q,k,l \leq r$ and $a,b \in \mfa$
					\begin{displaymath}
						[ \EUL{p}{q}(a), \ \EUL{k}{l}(b) ] = \delta_{q,k} \EUL{p}{l}(ab) + \delta_{l,p} \EUL{k}{q}(-ba).
					\end{displaymath}

		\item[]\label{ul ur}
			$[\ulbox, \ \urbox\phantom{}]$
				Given any $1 \leq p,q,k,l \leq r$ and $a,b \in \mfa$
					\begin{displaymath}
						[ \EUL{p}{q}(a), \ \EUR{k}{l}(b) ] = \delta_{q,k} \EUR{p}{l}(ab) + \delta_{q,l} \EUR{p}{k}(-a\ovlb\,).
					\end{displaymath}

		\item[]\label{ul bl}
			$[\ulbox, \ \blbox]$
				Given any $1 \leq p,q,k,l \leq r$ and $a,b \in \mfa$
					\begin{displaymath}
						[ \EUL{p}{q}(a), \ \EBL{k}{l}(b) ] = \delta_{p,k} \EBL{q}{l}(-\ovla\,b) + \delta_{p,l} \EBL{q}{k}(\ovlab\,).
					\end{displaymath}

		\item[]\label{ur ur}
			$[\urbox\phantom{}, \ \urbox\phantom{}]$
				Given any $1 \leq p,q,k,l \leq r$ and $a,b \in \mfa$
					\begin{displaymath}
						[ \EUR{p}{q}(a), \ \EUR{k}{l}(b) ] = 0.
					\end{displaymath}

		\item[]\label{ur bl}
			$[\urbox\phantom{}, \ \blbox]$
				Given any $1 \leq p,q,k,l \leq r$ and $a,b \in \mfa$
					\begin{align*}
						&[ \EUR{p}{q}(a), \ \EBL{k}{l}(b) ] = \\
						& \quad \delta_{p,k} \EUL{q}{l}(-\ovla\,b) + \delta_{p,l} \EUL{q}{k}(\ovlab\,) + \delta_{q,k} \EUL{p}{l}(ab) + \delta_{q,l} 	\EUL{p}{k}(-a\ovlb\,).
					\end{align*}

		\item[]\label{bl bl}
			$[\blbox, \ \blbox]$
				Given any $1 \leq p,q,k,l \leq r$ and $a,b \in \mfa$
					\begin{displaymath}
						[ \EBL{p}{q}(a), \ \EBL{k}{l}(b) ] = 0.
					\end{displaymath}

	\end{enumerate}
	\end{lemma*}

		\begin{propn}
			The module $C$ is zero.
		\end{propn}
		\begin{proof}
			The lemmas $ \left[ \vertbox, \ \vertbox \right], \quad \left[ \vertbox, \ \hortbox \right], \ldots, \left[ \blbox, \ \blbox \right] $ reveal that the image $\psi\left(\imAd\right)$ has trivial intersection with
				\begin{align*}
					& \left\{
							\left[
								\begin{array}{cc}
									0 & \chi_{e_i c} \\
									\left(e_{n+1-i} c\right)^t & 0
								\end{array}
							\right]
							\, : \ 
							c \in C
						\right\} \\
					& \quad \bigcup
						\left\{
							\left[
								\begin{array}{cc}
									0 & \chi_{e_{n+1-i} c} \\
									\left(e_i c\right)^t & 0
								\end{array}
							\right]
							\, : \ 
							c \in C
						\right\}
						\bigcup
						\left\{
							\left[
								\begin{array}{cc}
									0 & \chi_{e_{r+1} c} \\
									\left(e_{r+1} c\right)^t & 0
								\end{array}
							\right]
							\, : \ 
							c \in C
						\right\}
				\end{align*}
			Since the module $C$ is only involved in one of these subsets in the construction of $\sotnpoaecc$, and since our homomorphism $\psi: \imAd \to \sotnpoaecc$ is surjective, the triviality of the intersection implies that
				\begin{displaymath}
					C = \left\{ 0 \right\}.
				\end{displaymath}
		\end{proof}


\section{Achieving  a ``sufficient'' understanding of $\mfa$, $\eta$, $C$, and $\chi$}\label{sect:sufficient understanding}

	In the previous section we used the homomorphism $\psi: \imAd \to \sotnpoaecc$, given by [ABnG]'s Recognition Theorem, to get a sense (i) of what the generators of $\mfa$ ought to be; (ii) of what the involution $\eta$ on $\mfa$ ought to be; (iii) of what the relations on the generators of $\mfa$ ought to be; and (iv) that $C = 0$ and $\chi =0$.  
	
	But a key question remains: Could it be that the analysis in the previous section was incomplete in that it failed to detect other fundamental relations among the members of $\mfa$?  
		
	We settle this question as follows:  
		\begin{enumerate}
			\item
				Take the $4$-tuple of associative algebra, involution, module, and hermitian form as we presently understand it. That is,
					
					\begin{enumerate}
						\item[(i)]
							Let
								\begin{enumerate}
									\item[]
										$\Omega =$ the set of all long roots of the form $\pm \big( \eipej{i}{i+1} \big) $ that we have adjoined,
									\item[]
										$\Theta =$ the set of all long roots in $\Delta_B$ which we have adjoined but that are not in $\Omega$,
								\end{enumerate}
							and
								\begin{itemize}
									\item[]
										\begin{displaymath}
											X_e = \bigcup_{\omega \in \Omega} \left\{ x_{\omega,1}, \ldots, x_{\omega,\Nomega} \right\}
										\end{displaymath}
									\item[]
										\begin{displaymath}
											X_f = \bigcup_{\omega \in \Omega} \left\{ x^{-1}_{\omega,1}, \ldots, x^{-1}_{\omega,\Nomega} \right\}
										\end{displaymath}			
									\item[]
										\begin{displaymath}
											Y_e = \bigcup_{\theta \in \Theta} \left\{ y_{\theta,1}, \ldots, y_{\theta,\Ntheta} \right\}
										\end{displaymath}
									\item[]
										\begin{displaymath}
											Y_f = \bigcup_{\theta \in \Theta} \left\{ y^{-1}_{\theta,1}, \ldots, y^{-1}_{\theta,\Ntheta} \right\}
										\end{displaymath}						
									\item[]
										\begin{displaymath}
											Z_e = \bigcup_{\theta \in \Theta} \left\{ z_{\theta,1}, \ldots, z_{\theta,\Ntheta} \right\}
										\end{displaymath}
									\item[]
										\begin{displaymath}
											Z_f = \bigcup_{\theta \in \Theta} \left\{ z^{-1}_{\theta,1}, \ldots, z^{-1}_{\theta,\Ntheta} \right\}
										\end{displaymath}
								\end{itemize}	
							denote collections of indeterminates indexed by the sets $\Omega$ and $\Theta$.  Let $\mfb$ be the unital associative $\bbC$-algebra generated by the indeterminates in
		\begin{displaymath}
			X_e \cup X_f \cup Y_e \cup Y_f \cup Z_e \cup Z_f,
		\end{displaymath}
subject to the relations
		\begin{itemize}
			\item[]
				$ y_{\eimej{p}{q},i} \  x_{\eipej{p}{q},j} = x_{\eipej{p}{q},j} \  z_{\eimej{p}{q},i}$
			\item[]
				$ y_{\eimej{p}{q},i} \  z_{\eipej{p}{q},j} = y_{\eipej{p}{q},j} \  z_{\eimej{p}{q},i}$
			\item[]
				$ z_{\eimej{p}{q},i} \  x_{\meimej{p}{q},k} = x_{\meimej{p}{q},k} \  y_{\eimej{p}{q},i}$
			\item[]
				$ z_{\eimej{p}{q},i} \  y_{\meimej{p}{q},k} = z_{\meimej{p}{q},k} \  y_{\eimej{p}{q},i}$			
		\end{itemize}	
where $i = 1, \ldots, N_{\eimej{p}{q}}$ for $\eimej{p}{q} \in \Theta$,  $j =  1 \ldots N_{\eipej{p}{q}}$ for $\eipej{p}{q} \in \Omega \cup \Theta$, and $k = 1, \ldots, N_{\meimej{p}{q}}$ for $\meimej{p}{q} \in \Omega \cup \Theta$. 

		\item[(ii)]
			Define an involution, which we also call $\eta$ and sometimes denote by $\bar{\cdot}$, on $\mfb$ by 
				\begin{itemize}
					\item[]
						$\eta\left( x_{\omega,i} \right) = x_{\omega,i}$, if $\omega \in \Omega$ and $1 \leq i \leq \Nomega$,
					\item[]
						$\eta\left( y_{\theta,i} \right) = z_{\theta,i}$, if $\theta \in \Theta$ and $1 \leq i \leq \Ntheta$,
					\item[]
						$\eta\left( z_{\theta,i} \right)= y_{\theta,i}$, if $\theta \in \Theta$ and $1 \leq i \leq \Ntheta$.
					\end{itemize}
		\item[(iii)]
			Let $C = 0$ be the trivial $\mfb$-module.
		\item[(iv)]
			Let $\chi=0$ be a hermitian form on $C$.
				\end{enumerate}

\noindent \emph{\bf Remarks:}
	\begin{itemize}
		\item[(a)]
			The indeterminates in $X_e \cup X_f \cup \cdots \cup Z_f$ are intended to capture the elements of the form $a$, $a'$, $b$, $b'$, $c$, and $c'$ of $\mfa$ that we studied in Section \ref{sect:minimal understanding}, which arose from the images of the map $\psi$.
		\item[(b)]
			In the relations listed above, our use of the indeterminates $x_{\eipej{p}{q},j}$ and $x_{\meimej{p}{q},j}$ signals that we are working with roots in $\Omega$ and, hence, $|p-q| = 1$ in this setting.  Likewise, our use of the indeterminates $y_{\eipej{p}{q},j}$, $z_{\eipej{p}{q},j}$, $y_{\meimej{p}{q},j}$, and $z_{\meimej{p}{q},j}$ signals that we are working with roots in $\Theta$ and $p$, $q$ such that $|p-q| \geq 2$.  	
	\end{itemize}

			\item
				Construct a map 
					\begin{displaymath}
						\varphi: \gimAd \to \sotrpobecc
					\end{displaymath}
				sending the generators
					\begin{displaymath}
						e_1, \ldots, e_r, \ \bigcup_{\omega \in \Omega} \left\{ e_{\omega,1}, \ldots, e_{\omega,N_{\omega}} \right\}, \ \bigcup_{\theta \in \Theta} \left\{ e_{\theta,1}, \ldots, e_{\theta,N_{\theta}} \right\}
					\end{displaymath}
					\begin{displaymath}
						f_1, \ldots, f_r, \ \bigcup_{\omega \in \Omega} \left\{ f_{\omega,1}, \ldots, f_{\omega,N_{\omega}} \right\}, \ \bigcup_{\theta \in \Theta} \left\{ f_{\theta,1}, \ldots, f_{\theta,N_{\theta}} \right\}
					\end{displaymath}
					\begin{displaymath}
						h_1, \ldots, h_r, \ \bigcup_{\omega \in \Omega} \left\{ h_{\omega,1}, \ldots, h_{\omega,N_{\omega}} \right\}, \ \bigcup_{\theta \in \Theta} \left\{ h_{\theta,1}, \ldots, h_{\theta,N_{\theta}} \right\}
					\end{displaymath}
	
				of $\gimAd$ to

					\begin{displaymath}
						\te_1, \ldots, \te_r, \ \bigcup_{\omega \in \Omega} \left\{ \te_{\omega,1}, \ldots, \te_{\omega,N_{\omega}} \right\}, \ \bigcup_{\theta \in \Theta} \left\{ \te_{\theta,1}, \ldots, \te_{\theta,N_{\theta}} \right\}
					\end{displaymath}
					\begin{displaymath}
						\tf_1, \ldots, \tf_r, \ \bigcup_{\omega \in \Omega} \left\{ \tf_{\omega,1}, \ldots, \tf_{\omega,N_{\omega}} \right\}, \ \bigcup_{\theta \in \Theta} \left\{ \tf_{\theta,1}, \ldots, \tf_{\theta,N_{\theta}} \right\}
					\end{displaymath}
					\begin{displaymath}
						\tilh_1, \ldots, \tilh_r, \ \bigcup_{\omega \in \Omega} \left\{ \tilh_{\omega,1}, \ldots, \tilh_{\omega,N_{\omega}} \right\}, \ \bigcup_{\theta \in \Theta} \left\{ \tilh_{\theta,1}, \ldots, \tilh_{\theta,N_{\theta}} \right\}
					\end{displaymath}

				respectively, where
					\begin{itemize}

						\item[]
							$\te_i := \EUL{i}{i+1}(1)$, \quad $ 1 \leq i \leq r-1$,

						\item[]
							$\te_r := \Evert{r}(\sqrt{2})$,
	
						\item[]
							$\te_{\omega,i} := \left\{ 
								\begin{array}{ll} 
									\EUR{p}{(p+1)}(x_{\omega,i}) & \text{if } \omega = \eipej{p}{p+1} \\ 
									\EBL{p}{p+1}(x_{\omega,i}) & \text{if } \omega = \meimej{p}{p+1}  
								\end{array} 
							\right.,$ \\
\vskip2pt	
							for $\omega \in \Omega$ and $1 \leq i \leq \Nomega$, and 

						\item[]
							$\te_{\theta,i} := \left\{ 
								\begin{array}{ll} 
									\EUL{p}{q}(y_{\theta,i}) & \text{if } \theta = \eimej{p}{q} \\ 
									\EUR{p}{q}(y_{\theta,i}) & \text{if } \theta = \eipej{p}{q} \\ 
									\EBL{p}{q}(y_{\theta,i}) & \text{if } \theta = \meimej{p}{q}  
								\end{array} 
							\right.,$ \\
\vskip2pt	
							for $\theta \in \Theta$ and $1 \leq i \leq \Ntheta$.
		
					\end{itemize}

\vskip2cm

					\begin{itemize}

						\item[]
							$\tf_i := \EUL{i+1}{i}(1)$, \quad $ 1 \leq i \leq r-1$,

						\item[]
							$\tf_r := \Ehort{r}(\sqrt{2})$,

						\item[]
							$\tf_{\omega,i} := \left\{ 
								\begin{array}{ll} 
									\EBL{(p+1)}{p}(x^{-1}_{\omega,i}) & \text{if } \omega = \eipej{p}{p+1} \\ 
									\EUR{p+1}{p}(x^{-1}_{\omega,i}) & \text{if } \omega = \meimej{p}{p+1}  
								\end{array} 
							\right.,$ \\
\vskip2pt	
							for $\omega \in \Omega$ and $1 \leq i \leq \Nomega$, and 

						\item[]
							$\tf_{\theta,i} := \left\{ 
								\begin{array}{ll} 
									\EUL{q}{p}(y^{-1}_{\theta,i}) & \text{if } \theta = \eimej{p}{q} \\ 
									\EBL{q}{p}(y^{-1}_{\theta,i}) & \text{if } \theta = \eipej{p}{q} \\ 
									\EUR{q}{p}(y^{-1}_{\theta,i}) & \text{if } \theta = \meimej{p}{q}  
								\end{array} 
							\right.,$ \\
\vskip2pt	
							for $\theta \in \Theta$ and $1 \leq i \leq \Ntheta$.
		
					\end{itemize}

\vskip2cm

					\begin{itemize}

						\item[]
							$\tilh_i := \EUL{i}{i}(1) + \EUL{i+1}{i+1}$, \quad $ 1 \leq i \leq r-1$,

						\item[]
							$\tilh_r := \EUL{r}{r}(2)$,

						\item[]
							$\tilh_{\omega,i} := \left\{ 
								\begin{array}{ll} 
									\EUL{p}{p}(1) + \EUL{p+1}{p+1}(1) & \text{if } \omega = \eipej{p}{p+1} \\ 
									\EUL{p}{p}(-1) + \EUL{p+1}{p+1}(-1) & \text{if } \omega = \meimej{p}{p+1}  
								\end{array} 
							\right.,$ \\
\vskip2pt	
							for $\omega \in \Omega$ and $1 \leq i \leq \Nomega$, and 

						\item[]
							$\tilh_{\theta,i} := \left\{ 
								\begin{array}{ll} 
									\EUL{p}{p}(1) + \EUL{q}{q}(-1) & \text{if } \theta = \eimej{p}{q} \\ 
									\EUL{p}{p}(1) + \EUL{q}{q}(1) & \text{if } \theta = \eipej{p}{q} \\ 
									\EUL{p}{p}(-1) + \EUL{q}{q}(-1) & \text{if } \theta = \meimej{p}{q}  
								\end{array} 
							\right.,$ \\
\vskip2pt	
							for $\theta \in \Theta$ and $1 \leq i \leq \Ntheta$.
		
					\end{itemize}

			\item
				Show that $\varphi$ is
					\begin{enumerate}
						\item
							a Lie algebra homomorphism ($\S$\ref{subsect:is homomorphism}),
						\item
							that is surjective ($\S$\ref{subsect:is surjective}), and
						\item
							graded ($\S$\ref{subsect:is graded}).
					\end{enumerate}
			\item
				Show that the radical $\mfr$ of $\gimAd$ lies in the kernel of this map $\varphi$ ($\S$\ref{subsect:is graded}), hence inducing a surjective, graded, Lie algebra homomorphism
					\begin{displaymath}
						\phi: \imAd \to \sotrpobecc.
					\end{displaymath}	
				\item
					And finally, we show that $\phi$ is a central map ($\S$\ref{subsect:is central}).
		\end{enumerate} 
	
	The last two items reveal that the relations on $\mfa$ that we arrived at in $\S$\ref{subsect:relations} are in fact all the ones we needed.  That is, we really do have the correct coordinate algebra $\mfa$ in hand.

\subsection{}\label{subsect:is homomorphism}

			\begin{thm}\label{thm:very long theorem}
				The map $\varphi: \gimAd \to \sotrpobecc$ is a Lie algebra homomorphism.
			\end{thm}
			
\begin{proof}


We show that the images in $\sotrpobecc$ of the generators of $\gimAd$, under the map $\varphi$, satisfy the relations (R1) - (R3) of Definition \ref{defn:gim algebra} with respect to the same $(r+d) \times (r+d)$ generalized intersection matrix $A^{[d]}$ as used in the construction of the algebra $\gimAd$.

While working with the various long roots in our proof, we use labels like $u$ or $v$ to denote the indeterminates $x_{\omega,i}$ or $y_{\theta,i}$.

The reason that we can substitute $u$ or $v$ for the actual indeterminates is that the result of taking a bracket like
	\begin{displaymath}
		\left[ \te_{\meimej{p}{q},i}, \ \te_{\meimej{k}{l},j}  \right] = \left[ \EBL{p}{q}(y_{\meimej{p}{q},i}), \ \EBL{k}{l}(y_{\meimej{k}{l},j}) \right]
	\end{displaymath}
depends primarily on the indices $p$, $q$, $k$, and $l$ rather than on the particular elements of the algebra $\mfb$ being housed at these sites.

If we agree on this convention of using substitute variables like $u$, then we must recognize that 
	\begin{displaymath}
		\ovlu = \left\{
		\begin{array}{ll}
		x_{\omega,i} & \text{ if } u = x_{\omega,i} \\
		z_{\theta,i} & \text{ if } u = y_{\theta,j} 
		\end{array}
		\right.
	\end{displaymath}
That is, the involution $\bar{\cdot}$ applied to $u$ depends on whether $u$ is substituting for a variable associated to a root in $\Omega$ or a root in $\Theta$.



We show the computations for the interactions between the generators corresponding to the long roots  $\epsilon_p - \epsilon_q$, and $\epsilon_k - \epsilon_l$.  The remaining computations are similar. 

\vskip5pt

Let $ 1 \leq p,q,k,l \leq r$ with $p \neq q$ and $k \neq l$;  $u,v \in \{x_{\omega,i},  \, x_{\omega,j},  \, y_{\theta,i} \, y_{\theta,j} \}$ and $\uinv, \vinv\in \{x^{-1}_{\omega,i}, \, x^{-1}_{\omega,j}, \, y^{-1}_{\theta,i} \, y^{-1}_{\theta,j} \}$, where $\omega \in \Omega$, $\theta \in \Theta$, and $1 \leq i,j \leq \Nomega$ or $1 \leq i,j \leq \Ntheta$.

Using the definition of $A_{\eimej{p}{q}, \eimej{k}{l}}$, we see that 
	\begin{displaymath}
		A_{\eimej{p}{q}, \eimej{k}{l}} = \delta_{p,k} - \delta_{p,l} - \delta_{q,k} + \delta_{q,l}		
	\end{displaymath}
	\begin{displaymath}
		= \left\{ 
		\begin{array}{rl}
			0 & \text{ if } p,q \notin \{k,l\} \\
			1 & \text{ if } p=k \text{ but } q \neq l \\
			-1 & \text{ if } p=l \text{ but } q \neq k \\
			-1 & \text{ if } p \neq l \text{ but } q = k \\
			1 & \text{ if } p \neq k \text{ but } q = l \\
			2 & \text{ if } p=k \text{ and } q = l \\
			-2 & \text{ if } p=l \text{ and } q = k \\
		\end{array}
		\right.
	\end{displaymath}


\noindent A.
	\begin{align*}
		& \left[ \te_{\eimej{p}{q},i}, \ \te_{\eimej{k}{l},j}   \right] = \left[ \EUL{p}{q}(u), \ \EUL{k}{l}(v) \right] \\
		& = 
			\delta_{q,k} \EUL{p}{l} \left( uv \right) +
			\delta_{l,p} \EUL{k}{q} \left( -vu \right) \\
	\end{align*}
	\begin{displaymath}
		= \left\{ 
			\begin{array}{ll}
				0 & \text{ if } p,q \notin \{k,l\} \\
				0 & \text{ if } p=k \text{ but } q \neq l \\
				\EUL{k}{q} \left( -vu \right) & \text{ if } p=l \text{ but } q \neq k \\
				\EUL{p}{l} \left( uv \right) & \text{ if } p \neq l \text{ but } q = k \\
				0 & \text{ if } p \neq k \text{ but } q = l \\
				0 & \text{ if } p=k \text{ and } q = l \\
				\EUL{p}{p} \left( uv \right) + \EUL{q}{q} \left( - vu \right) & \text{ if } p=l \text{ and } q = k \\
			\end{array}
		\right.
	\end{displaymath}
	
	\begin{itemize}
		\item
			If $p=l$ but $q \neq k$, then $ \left[ \EUL{p}{q}(u), \ \EUL{k}{q} \left( - vu \right) \right] = 0$ because $q \neq k$ and $q \neq p$.
		\item
			If $p \neq l$ but $q = k$, then $ \left[ \EUL{p}{q}(u), \ \EUL{p}{l} \left( uv \right) \right] = 0$ because $q \neq p$ and $l \neq p$.
		\item
			If $p = l$ and $q = k$, then 
				\begin{align*}
					& \left[ \EUL{p}{q}(u), \ \EUL{p}{p} \left( uv \right) + \EUL{q}{q} \left( -vu \right) \right] \\
					& = \EUL{p}{q}\left( - uvu \right)	+ \EUL{p}{q} \left( - uvu \right) \\			
					& = \EUL{p}{q} \left( -2 uvu \right).
				\end{align*}			
	\end{itemize}

So
	\begin{displaymath}
		\Big( \ad \te_{\eimej{p}{q},i} \Big)^{1+1} \ \te_{\eimej{k}{l},j}	= \left\{ 
			\begin{array}{ll}
				0 & \text{ if } p,q \notin \{k,l\} \\
				0 & \text{ if } p=k \text{ but } q \neq l \\
				0 & \text{ if } p=l \text{ but } q \neq k \\
				0 & \text{ if } p \neq l \text{ but } q = k \\
				0 & \text{ if } p \neq k \text{ but } q = l \\
				0 & \text{ if } p=k \text{ and } q = l \\
				\EUL{p}{q} \left( -2 uvu \right) & \text{ if } p=l \text{ and } q = k \\
			\end{array}
		\right.
	\end{displaymath}

Since $\left[ \EUL{p}{q}(u), \ \EUL{p}{q} \left( -2 uvu \right) \right] = 0$, we get that
	\begin{displaymath}
		\Big( \ad \te_{\eimej{p}{q},i} \Big)^{2+1} \ \te_{\eimej{k}{l},j} = 0.
	\end{displaymath}


\noindent B.
	\begin{align*}
		& \left[ \tf_{\eimej{p}{q},i}, \ \tf_{\eimej{k}{l},j}   \right] = \left[ \EUL{q}{p}(\uinv), \ \EUL{l}{k}(\vinv) \right] \\
		& = 
			\delta_{p,l} \EUL{q}{k} \left( \uinv \vinv \right) +
			\delta_{k,q} \EUL{l}{p} \left( -\vinv \uinv \right) \\
	\end{align*}
	\begin{displaymath}
		= \left\{ 
			\begin{array}{ll}
				0 & \text{ if } p,q \notin \{k,l\} \\
				0 & \text{ if } p=k \text{ but } q \neq l \\
				\EUL{q}{k} \left( \uinv \vinv \right) & \text{ if } p=l \text{ but } q \neq k \\
				\EUL{l}{p} \left( -\vinv \uinv \right) & \text{ if } p \neq l \text{ but } q = k \\
				0 & \text{ if } p \neq k \text{ but } q = l \\
				0 & \text{ if } p=k \text{ and } q = l \\
				\EUL{p}{p} \left( -\vinv\uinv \right) + \EUL{q}{q} \left( \uinv \vinv \right) & \text{ if } p=l \text{ and } q = k \\
			\end{array}
		\right.
	\end{displaymath}
	
	\begin{itemize}
		\item
			If $p=l$ but $q \neq k$, then $ \left[ \EUL{q}{p}(\uinv), \ \EUL{q}{k} \left( \uinv \vinv \right) \right] = 0$ because $p \neq q$ and $k \neq q$.
		\item
			If $p \neq l$ but $q = k$, then $ \left[ \EUL{q}{p}(\uinv), \ \EUL{l}{p} \left( -\vinv \uinv \right) \right] = 0$ because $p \neq l$ and $p \neq q$.
		\item
			If $p = l$ and $q = k$, then 
				\begin{align*}
					& \left[ \EUL{q}{p}(\uinv), \ \EUL{p}{p} \left( -\vinv \uinv \right) + \EUL{q}{q} \left( \uinv \vinv \right) \right] \\
					& = \EUL{q}{p}\left( - \uinv \vinv \uinv \right)	+ \EUL{q}{p} \left( - \uinv \vinv \uinv \right) \\			
					& = \EUL{q}{p} \left( -2 \uinv \vinv \uinv \right).
				\end{align*}			
	\end{itemize}

So
	\begin{displaymath}
		\Big( \ad \tf_{\eimej{p}{q},i} \Big)^{1+1} \ \tf_{\eimej{k}{l},j}	= \left\{ 
			\begin{array}{ll}
				0 & \text{ if } p,q \notin \{k,l\} \\
				0 & \text{ if } p=k \text{ but } q \neq l \\
				0 & \text{ if } p=l \text{ but } q \neq k \\
				0 & \text{ if } p \neq l \text{ but } q = k \\
				0 & \text{ if } p \neq k \text{ but } q = l \\
				0 & \text{ if } p=k \text{ and } q = l \\
				\EUL{q}{p} \left( -2 \uinv \vinv \uinv \right) & \text{ if } p=l \text{ and } q = k \\
			\end{array}
		\right.
	\end{displaymath}

Since $\left[ \EUL{q}{p}(\uinv), \ \EUL{q}{p} \left( -2 \uinv \vinv \uinv \right) \right] = 0$, we get that
	\begin{displaymath}
		\Big( \ad \tf_{\eimej{p}{q},i} \Big)^{2+1} \ \tf_{\eimej{k}{l},j} = 0.
	\end{displaymath}


\noindent C.
	\begin{align*}
		& \left[ \tilh_{\eimej{p}{q},i}, \ \tilh_{\eimej{k}{l},j}   \right] \\
		& = \left[ \EUL{p}{p}(1) + \EUL{q}{q}(-1), \ \EUL{k}{k}(1) + \EUL{l}{l}(-1) \right] \\
		& = 
			\delta_{p,k} \EUL{p}{p} \left( [1,1] \right) +
			\delta_{p,l} \EUL{p}{p} \left( [1,-1] \right) +
			\delta_{q,k} \EUL{q}{q} \left( [-1,1] \right) +
			\delta_{q,l} \EUL{q}{q} \left( [-1,-1] \right) \\
		& = 
			\delta_{p,k} \EUL{p}{p} \left( 0 \right) +
			\delta_{p,l} \EUL{p}{p} \left( 0 \right) +
			\delta_{q,k} \EUL{q}{q} \left( 0 \right) +
			\delta_{q,l} \EUL{q}{q} \left( 0 \right) \\
		& = 0
	\end{align*}


\noindent D.
	\begin{align*}
		& \left[ \te_{\eimej{p}{q},i}, \ \tf_{\eimej{k}{l},j}   \right] = \left[ \EUL{p}{q}(u), \ \EUL{l}{k}(\vinv) \right] \\
		& = 
			\delta_{q,l} \EUL{p}{k} \left( u \, \vinv \right) +
			\delta_{k,p} \EUL{l}{q} \left( -\vinv \, u \right) \\
	\end{align*}
	\begin{displaymath}
		= \left\{ 
			\begin{array}{ll}
				0 & \text{ if } p,q \notin \{k,l\} \\
				\EUL{l}{q} \left( -\vinv \, u \right) & \text{ if } p=k \text{ but } q \neq l \\
				0 & \text{ if } p=l \text{ but } q \neq k \\
				0 & \text{ if } p \neq l \text{ but } q = k \\
				\EUL{p}{k} \left( u \, \vinv \right) & \text{ if } p \neq k \text{ but } q = l \\
				\EUL{p}{p} \left( u \, \vinv \right) + \EUL{q}{q} \left( -\vinv \, u \right) & \text{ if } p=k \text{ and } q = l \\
				0 & \text{ if } p=l \text{ and } q = k \\
			\end{array}
		\right.
	\end{displaymath}
	
	\begin{itemize}
		\item
			If $p=k$ but $q \neq l$, then $ \left[ \EUL{p}{q}(u), \ \EUL{l}{q} \left( -\vinv \, u \right) \right] = 0$ because $q \neq l$ and $q \neq p$.
		\item
			If $p \neq k$ but $q = l$, then $ \left[ \EUL{p}{q}(u), \ \EUL{p}{k} \left( u \, \vinv \right) \right] = 0$ because $q \neq p$ and $k \neq p$.
		\item
			If $p = k$ and $q = l$, then 
				\begin{align*}
					& \left[ \EUL{p}{q}(u), \ \EUL{p}{p} \left( u \, \vinv \right) + \EUL{q}{q} \left( -\vinv u \right) \right] \\
					& = \EUL{p}{q}\left( - u \, \vinv \, u \right)	+ \EUL{p}{q} \left( - u \, \vinv \, u \right) \\			
					& = \EUL{p}{q} \left( -2 u \, \vinv \, u \right).
				\end{align*}			
	\end{itemize}

So
	\begin{displaymath}
		\Big( \ad \te_{\eimej{p}{q},i} \Big)^{1+1} \ \tf_{\eimej{k}{l},j}	= \left\{ 
			\begin{array}{ll}
				0 & \text{ if } p,q \notin \{k,l\} \\
				0 & \text{ if } p=k \text{ but } q \neq l \\
				0 & \text{ if } p=l \text{ but } q \neq k \\
				0 & \text{ if } p \neq l \text{ but } q = k \\
				0 & \text{ if } p \neq k \text{ but } q = l \\
				\EUL{p}{q} \left( -2 u \, \vinv \, u \right) & \text{ if } p=k \text{ and } q = l \\
				0 & \text{ if } p=l \text{ and } q = k \\
			\end{array}
		\right.
	\end{displaymath}

Since $\left[ \EUL{p}{q}(u), \ \EUL{p}{q} \left( -2 u \, \vinv \, u \right) \right] = 0$, we get that
	\begin{displaymath}
		\Big( \ad \te_{\eimej{p}{q},i} \Big)^{2+1} \ \tf_{\eimej{k}{l},j} = 0.
	\end{displaymath}


\noindent E.
	\begin{align*}
		& \left[ \tf_{\eimej{p}{q},i}, \ \te_{\eimej{k}{l},j}   \right] = \left[ \EUL{q}{p}(\uinv), \ \EUL{k}{l}(v) \right] \\
		& = 
			\delta_{p,k} \EUL{q}{l} \left( \uinv v \right) +
			\delta_{l,q} \EUL{k}{p} \left( -v \uinv \right) \\
	\end{align*}
	\begin{displaymath}
		= \left\{ 
			\begin{array}{ll}
				0 & \text{ if } p,q \notin \{k,l\} \\
				\EUL{q}{l} \left( \uinv \, v \right) & \text{ if } p=k \text{ but } q \neq l \\
				0 & \text{ if } p=l \text{ but } q \neq k \\
				0 & \text{ if } p \neq l \text{ but } q = k \\
				\EUL{k}{p} \left( -v \, \uinv \right) & \text{ if } p \neq k \text{ but } q = l \\
				\EUL{p}{p} \left( -v \, \uinv \right) + \EUL{q}{q} \left( \uinv \, v \right) & \text{ if } p=k \text{ and } q = l \\
				0 & \text{ if } p=l \text{ and } q = k \\
			\end{array}
		\right.
	\end{displaymath}
	
	\begin{itemize}
		\item
			If $p=k$ but $q \neq l$, then $ \left[ \EUL{q}{p}(\uinv), \ \EUL{q}{l} \left( \uinv \, v \right) \right] = 0$ because $p \neq q$ and $l \neq q$.
		\item
			If $p \neq k$ but $q = l$, then $ \left[ \EUL{q}{p}(\uinv), \ \EUL{k}{p} \left( -v \, \uinv \right) \right] = 0$ because $p \neq k$ and $p \neq q$.
		\item
			If $p = k$ and $q = l$, then 
				\begin{align*}
					& \left[ \EUL{q}{p}(\uinv), \ \EUL{p}{p} \left( -v \, \uinv \right) + \EUL{q}{q} \left( \uinv \, v \right) \right] \\
					& = \EUL{q}{p}\left( - \uinv \, v \, \uinv \right)	+ \EUL{q}{p} \left( - \uinv \, v \, \uinv \right) \\			
					& = \EUL{q}{p} \left( -2 \uinv \, v \, \uinv \right).
				\end{align*}			
	\end{itemize}

So
	\begin{displaymath}
		\Big( \ad \tf_{\eimej{p}{q},i} \Big)^{1+1} \ \te_{\eimej{k}{l},j}	= \left\{ 
			\begin{array}{ll}
				0 & \text{ if } p,q \notin \{k,l\} \\
				0 & \text{ if } p=k \text{ but } q \neq l \\
				0 & \text{ if } p=l \text{ but } q \neq k \\
				0 & \text{ if } p \neq l \text{ but } q = k \\
				0 & \text{ if } p \neq k \text{ but } q = l \\
				\EUL{q}{p} \left( -2 \uinv \, v \, \uinv \right) & \text{ if } p=k \text{ and } q = l \\
				0 & \text{ if } p=l \text{ and } q = k \\
			\end{array}
		\right.
	\end{displaymath}

Since $\left[ \EUL{q}{p}(\uinv), \ \EUL{q}{p} \left( -2 \uinv \, v \, \uinv \right) \right] = 0$, we get that
	\begin{displaymath}
		\Big( \ad \tf_{\eimej{p}{q},i} \Big)^{2+1} \ \te_{\eimej{k}{l},j} = 0.
	\end{displaymath}


\noindent F.
	\begin{align*}
		& \left[ \tilh_{\eimej{p}{q},i}, \ \te_{\eimej{k}{l},j}   \right] \\
		& = \left[ \EUL{p}{p}(1) + \EUL{q}{q}(-1), \ \EUL{k}{l}(v) \right] \\
		& = 
			\delta_{p,k} \EUL{p}{l}(v) +
			\delta_{l,p} \EUL{k}{p}(-v) +
			\delta_{q,k} \EUL{q}{l}(-v) +
			\delta_{l,q} \EUL{k}{q}(v) \\
		& = 
			\delta_{p,k} \EUL{k}{l}(v) -
			\delta_{p,l} \EUL{k}{l}(v) -
			\delta_{q,k} \EUL{k}{l}(v) +
			\delta_{q,l} \EUL{k}{l}(v) \\			
		& = \Big( \delta_{p,k} - \delta_{p,l} - \delta_{q,k} + \delta_{q,l} \Big) \EUL{k}{l}(v) \\
		& = A_{\eimej{p}{q}, \eimej{k}{l}} \ \  \te_{\eimej{k}{l},j}
	\end{align*}


\noindent G.
	\begin{align*}
		& \left[ \tilh_{\eimej{p}{q},i}, \ \tf_{\eimej{k}{l},j}   \right] \\
		& = \left[ \EUL{p}{p}(1) + \EUL{q}{q}(-1), \ \EUL{l}{k}\left( \vinv \right) \right] \\
		& = 
			\delta_{p,l} \EUL{p}{k}\left( \vinv \right) +
			\delta_{k,p} \EUL{l}{p}\left( -\vinv \right) +
			\delta_{q,l} \EUL{q}{k}\left( -\vinv \right) +
			\delta_{k,q} \EUL{l}{q}\left( \vinv \right) \\
		& = 
			\delta_{p,l} \EUL{l}{k}\left( \vinv \right) -
			\delta_{p,k} \EUL{l}{k}\left( \vinv \right) -
			\delta_{q,l} \EUL{l}{k}\left( \vinv \right) +
			\delta_{q,k} \EUL{l}{k}\left( \vinv \right) \\			
		& = - \Big( \delta_{p,k} - \delta_{p,l} - \delta_{q,k} + \delta_{q,l} \Big) \EUL{l}{k}\left( \vinv \right) \\
		& = - A_{\eimej{p}{q}, \eimej{k}{l}} \  \  \tf_{\eimej{k}{l},j}
	\end{align*}


\end{proof}



\subsection{}\label{subsect:is surjective}
	
	Our main goal in this subsection is to demonstrate that $\varphi: \gimAd \to \sotrpobecc$, as given above, is surjective.  
	
	Let $\mfB$ be the Lie subalgebra of $\sotrpobecc$ generated by the images, under $\varphi$, of the generators of $\gimAd$.  We assume that $d \geq 1$.  
	
	An arbitrary element in $\mfb$ is a finite $\bbC$-linear combination of monomials of the form 
		\begin{displaymath}
			\Tmanyl,
		\end{displaymath}
	where $l \geq 1$ and for $1 \leq j \leq l$, 
		\begin{itemize}
			\item[]
				$\mu_j \in \Omega \cup \Theta$,
			\item[]
				$k_j$ refers to the $k_j^{\text{th}}$ copy of the long root $\mu_j$ that we have adjoined, with $1 \leq k_j \leq N_{\mu_j}$,
			\item[]
				$m_j \in \bbZ$, and 
			\item[]
				$t_{\mu_j,k_j}$ is one of the indeterminates $x_{\mu_j,k_j}$, $y_{\mu_j,k_j}$, or $z_{\mu_j,k_j}$.
		\end{itemize}

	Since $\mfB$ is a vector space over $\bbC$, it suffices to demonstrate that given such a monomial and any $1 \leq i, j \leq r$,
		\begin{itemize}
			\item[]
				$\Evert{i} \left( \Tmanyl \right) \in \mfB$,
			\item[]
				$\Ehort{i} \left( \Tmanyl \right) \in \mfB$,
			\item[]
				$\EUL{i}{j} \left( \Tmanyl \right) \in \mfB$,								
			\item[]
				$\EUR{i}{j} \left( \Tmanyl \right) \in \mfB$, and				
			\item[]
				$\EBL{i}{j} \left( \Tmanyl \right) \in \mfB$.				
		\end{itemize}

	To show this, we need to do some preliminary work.
	
	\begin{lemma}\label{lemma:vertical a's}
		\begin{enumerate}
			\item
				If $\Evert{i}(a) \in \mfB$ for some $2 \leq i \leq r$, $a \in \mfb$, then $\Evert{i-1}(a) \in \mfB$.
			\item	
				If $\Evert{i}(a) \in \mfB$ for some $1 \leq i \leq r-1$, $a \in \mfb$, then $\Evert{i+1}(a) \in \mfB$.
		\end{enumerate}
	\end{lemma}
	
	\begin{proof}
		\begin{enumerate}
			\item
				We know that $\varphi\left( e_{i-1} \right) = \te_{i-1} = 		\EUL{i-1}{i}(1) \in \mfB$.  So
					\begin{align*}
						& [\EUL{i-1}{i}(1), \ \Evert{i}(a) ] \\
						& \ = \Evert{i-1}(a)
					\end{align*}
				is also an element of $\mfB$.
			\item
					Since $\varphi\left( f_i \right) = \tf_{i} = 		\EUL{i+1}{i}(1)\in \mfB$, so is 
					\begin{displaymath}
						\Evert{i+1}(a) = \left[ \EUL{i+1}{i}(1), \ \Evert{i}(a) \right].
					\end{displaymath}
		\end{enumerate}		
	\end{proof}

	\begin{cor}\label{cor:vertical a's}
		If $ \Evert{i}(a) \in \mfB$ for some $1 \leq i \leq r$ and $a \in \mfb$, then $\Evert{j}(a) \in \mfB$ for all $j = 1, \ldots, r$. 
	\end{cor}

	\begin{cor}\label{cor:vertical ones}
		The element $\Evert{i}(1)$ lies in $\mfB$ for all $1 \leq i \leq r$.
	\end{cor}

	\begin{proof}
		$\varphi\left( e_r \right) = \te_r = \Evert{r}\left( \sqrt{2} \right)$ is an element of $\mfB$.  Hence, after multiplying by $\frac{1}{\sqrt{2}}$, $\Evert{r}(1) \in \mfB$.  By Lemma \ref{lemma:vertical a's}.1, it follows that $\Evert{i}(1) \in \mfB$ for all $1 \leq i < r$.
	\end{proof}

	\begin{lemma}\label{lemma:horizontal a's}
		\begin{enumerate}
			\item
				If $\Ehort{i}(a) \in \mfB$ for some $2 \leq i \leq r$, $a \in \mfb$, then $\Ehort{i-1}(a) \in \mfB$.
			\item	
				If $\Ehort{i}(a) \in \mfB$ for some $1 \leq i \leq r-1$, $a \in \mfb$, then $\Ehort{i+1}(a) \in \mfB$.
		\end{enumerate}
	\end{lemma}

	\begin{proof}
		\begin{enumerate}
			\item
				Since $\varphi\left( f_{i-1} \right) = \tf_{i-1} = 	 \EUL{i}{i-1}(1) \in \mfB$, so is
					\begin{displaymath}
						\Ehort{i-1}(a) = \left[ \Ehort{i}(a), \ \EUL{i}{i-1}(1) \right].
					\end{displaymath}
			\item
				Since $\varphi\left( e_i \right) = \te_{i} =	\EUL{i}{i+1}(1)  \in \mfB$, so is 
					\begin{displaymath}
						\Ehort{i+1}(a) = \left[ \Ehort{i}(a), \ \EUL{i}{i+1}(1) \right].
					\end{displaymath}
		\end{enumerate}		
	\end{proof}

	\begin{cor}\label{cor:horizontal a's}
		If $\Ehort{i}(a) \in \mfB$ for some $1 \leq i \leq r$, $a \in \mfb$, then $\Ehort{j}(a) \in \mfB$ for all $j = 1, \ldots, r$.
	\end{cor}

	\begin{cor}\label{cor:horizontal ones} 
		The element $\Ehort{i}(1)$ lies in $\mfB$ for all $1 \leq i \leq r$.
	\end{cor}

	\begin{cor}{\label{cor:diagonal ones}} 
		The element $\EUL{i}{i}(1)$ lies in $\mfB$ for all $1 \leq i \leq r$.
	\end{cor}

	\begin{proof}
		Take any $i = 1, \ldots, r$.  By Corollary \ref{cor:vertical ones}, $\Evert{i}(1) \in \mfB$.  By Corollary \ref{cor:horizontal ones}, $\Ehort{i}(1) \in \mfB$.  But then
		\begin{align*}
			& [ \Evert{i}(1), \ \Ehort{i}(1)  ] \\
			& \ = \EUL{i}{i}(1) \in \mfB.
		\end{align*}
	\end{proof}
	
		\begin{lemma}\label{lemma:eia to hia and meia} 
		If $\Evert{i}(a) \in \mfB$ for some $i = 1, \ldots, r$ and $a \in \mfb$, then $\EUL{j}{j}(a) \in \mfB$ and $\Ehort{j}(a) \in \mfB$ for all $j = 1, \ldots, r$.
	\end{lemma}

	\begin{proof}
		If $\Evert{i}(a) \in \mfB$ for some $i \in \{1, \ldots, r\}$ then, by Corollary \ref{cor:vertical a's}, $\Evert{j}(a) \in \mfB$ for all $j = 1, \ldots, r$.  We also, by Corollary \ref{cor:horizontal ones}, know that $\Ehort{j}(1) \in \mfB$ for all $j = 1, \ldots, r$.  But then, given any $j \in \{1, \ldots, r\}$, 
			\begin{displaymath}
				\left[ \Evert{j}(a), \ \Ehort{j}(1) \right] = \EUL{j}{j}(a) \in \mfB.
			\end{displaymath}
		This, in turn, implies that
			\begin{displaymath}
				\left[ \EUL{j}{j}(a), \ \Ehort{j}(1) \right] = \Ehort{j}(a) \in \mfB.
			\end{displaymath}
	\end{proof}

	\begin{cor}\label{cor:eia to an}
	If $\Evert{i}(a) \in \mfB$ for some $i = 1, \ldots, r$ and $a \in \mfb$, then $\Evert{j}\left( a^n \right)$, $\EUL{j}{j}\left( a^n \right)$, and $\Ehort{j}\left( a^n \right) \in \mfB$ for every positive integer $n$ and every $j = 1, \ldots, r$.
	\end{cor}
	
	\begin{proof}
		Lemma \ref{lemma:eia to hia and meia} tells us that it suffices to show $\Evert{i}\left( a^n \right) \in \mfB$ for every positive integer $n$.  We proceed by induction on $n$.  The base case holds because, by hypothesis, $\Evert{i}(a) \in \mfB$.  Suppose that $\Evert{i}\left( a^n \right) \in \mfB$ for some $n \geq 1$.  The hypothesis along with Lemma \ref{lemma:eia to hia and meia} also tells us that $\EUL{i}{i}(a) \in \mfB$.  But then
			\begin{displaymath}
				\left[ \EUL{i}{i}(a), \ \Evert{i}\left( a^n \right) \right] = \Evert{i}\left( a^{n+1} \right) \in \mfB.
			\end{displaymath}
	\end{proof}

	\begin{lemma}\label{lemma:meia to hia and eia}
		If $\Ehort{i}(a) \in \mfB$ for some $i = 1, \ldots, r$ and $a \in \mfb$, then $\EUL{j}{j}(a) \in \mfB$ and $\Evert{j}(a) \in \mfB$ for all $j = 1, \ldots, r$.
	\end{lemma}

	\begin{proof}
		If $\Ehort{i}(a) \in \mfB$ for some $i \in \{1, \ldots, r\}$ then, by Corollary \ref{cor:horizontal a's}, $\Ehort{j}(a) \in \mfB$ for all $j = 1, \ldots, r$.  We also, by Corollary \ref{cor:vertical ones}, know that $\Evert{j}(1) \in \mfB$ for all $j = 1, \ldots, r$.  But then, given any $j \in \{1, \ldots, r\}$, 
			\begin{displaymath}
				\left[ \Evert{j}(1), \ \Ehort{j}(a) \right] = \EUL{j}{j}(a) \in \mfB.
			\end{displaymath}
		This, in turn, implies that
			\begin{displaymath}
				\left[ \EUL{j}{j}(a), \ \Evert{j}(1) \right] = \Evert{j}(a) \in \mfB.
			\end{displaymath}
	\end{proof}

	\begin{cor}\label{cor:meia to an} 
		If $\Ehort{i}(a) \in \mfB$ for some $i = 1, \ldots, r$ and $a \in \mfb$, then $\Ehort{j}\left( a^n \right)$, $\EUL{j}{j}\left( a^n \right)$, and $\Evert{j}\left( a^n \right)$ all lie in $\mfB$ for every positive integer $n$ and every $j = 1, \ldots, r$.
	\end{cor}
	
	\begin{lemma}\label{lemma:off centre quads} 
		\begin{enumerate}
			\item
				If $\Evert{i}(a) \in \mfB$ and $\Ehort{j}(b) \in \mfB$ for some $1 \leq i,j \leq r$ with $i \neq j$ and $a, \, b \in \mfb$, then $\EUL{i}{j}(ab) \in \mfB$.
			\item
				If $\Evert{i}(a) \in \mfB$ and $\Evert{j}(b) \in \mfB$ for some $1 \leq i,j \leq r$ and $a, \, b \in \mfb$, then $\EUR{i}{j}(ab) \in \mfB$.
			\item
				If $\Ehort{i}(a) \in \mfB$ and $\Ehort{j}(b) \in \mfB$ for some $1 \leq i,j \leq r$ and $a, \, b \in \mfb$, then $\EBL{i}{j}(ab) \in \mfB$.
		\end{enumerate}
	\end{lemma}

	\begin{proof}
		\begin{enumerate}
			\item
				Since $\Evert{i}(a) \in \mfB$ and $\Ehort{j}(b) \in \mfB$,
					\begin{displaymath}
						\left[ \Evert{i}(a), \ \Ehort{j}(b) \right] = \EUL{i}{j}(ab) \in \mfB.			
					\end{displaymath}
			
			\item
				Since $\Evert{i}(a) \in \mfB$ and $\Evert{j}(b) \in \mfB$,
					\begin{displaymath}
						\left[ \Evert{j}(b), \ \Evert{i}(a) \right] = \EUR{i}{j}(ab) \in \mfB.			
					\end{displaymath}
			
			\item
				Since $\Ehort{i}(a) \in \mfB$ and $\Ehort{j}(b) \in \mfB$,
					\begin{displaymath}
						\left[ \Ehort{j}(b), \ \Ehort{i}(a) \right] = \EBL{i}{j}(ab) \in \mfB.			
					\end{displaymath}
		\end{enumerate}
	\end{proof}
	
	\begin{cor}\label{cor:off centre ones}
		The elements $\EUL{i}{j}(1)$, $\EUR{i}{j}(1)$, and $\EBL{i}{j}(1)$ lie in $\mfB$ for all $1 \leq i,j \leq r$.  
	\end{cor}
	
	Another point we must address before proceeding further is how we can get monomials involving indeterminates in $Z_e \cup Z_f$ given that, at least superficially, $\varphi$ sends the generators of $\gimAd$ to images involving indeterminates only in $X_e \cup X_f \cup Y_e \cup Y_f$.  
	
	Indeed,
		\begin{enumerate}
			\item
				suppose $\EUL{p}{q}\left( y_{\eimej{p}{q},i} \right) \in \mfB$, where $i$ refers to the $i^{\text{th}}$ copy of the long root $\eimej{p}{q}$ that we have adjoined, and $1 \leq p, q \leq r$ with $p \neq q$.  Then, having deduced that $\Evert{q}(-1)$, $\EBL{p}{q}(1)$, and $\Evert{p}(1) \in \mfB$, we get that
					\begin{align*}
						& \left[ \Evert{p}(1), \ \left[ \EBL{p}{q}(1), \ \left[ \Evert{q}(-1), \ \EUL{p}{q}\left( y_{\eimej{p}{q},i} \right) \right] \right] \right] \\
						& = \left[ \Evert{p}(1), \ \left[ \EBL{p}{q}(1), \ \Evert{p}\left( y_{\eimej{p}{q},i} \right) \right] \right] \\
						& = \left[ \Evert{p}(1), \ \Ehort{q}\left( \ovly_{\eimej{p}{q},i} \right) \right] \\
						& = \EUL{p}{q}\left( \ovly_{\eimej{p}{q},i} \right) \\
						& = \EUL{p}{q}\left( z_{\eimej{p}{q},i} \right) \in \mfB. \\
					\end{align*}
					
			\item
				suppose $\EUR{p}{q}\left( y_{\eipej{p}{q},i} \right) \in \mfB$, where $i$ refers to the $i^{\text{th}}$ copy of the long root $\eipej{p}{q}$ that we have adjoined, and $1 \leq p, q \leq r$ with $p \neq q$.  Then
					\begin{align*}
						& \left[ \Evert{q}(1), \ \left[ \Evert{q}(1), \ \left[ \Evert{p}(1), \ \left[ \EBL{p}{q}(1), \ \left[ \Ehort{q}(-1), \ \EUR{p}{q}\left( y_{\eipej{p}{q},i} \right) \right] \right] \right] \right] \right] \\
						& = \EUR{q}{p}\left( - y_{\eipej{p}{q},i} \right) \\
						& = \EUR{p}{q}\left( \ovly_{\eipej{p}{q},i} \right) \\
						& = \EUR{p}{q}\left( z_{\eipej{p}{q},i} \right) \in \mfB. \\
					\end{align*}					
			
						\item
				suppose $\EBL{p}{q}\left( y_{\meimej{p}{q},i} \right) \in \mfB$, where $i$ refers to the $i^{\text{th}}$ copy of the long root $\meimej{p}{q}$ that we have adjoined, and $1 \leq p, q \leq r$ with $p \neq q$.  Then
					\begin{align*}
						& \left[ \Ehort{p}(1), \ \left[ \Ehort{p}(1), \ \left[ \Ehort{q}(1), \ \left[ \EUR{p}{q}(1), \ \left[ \Evert{p}(1), \ \EBL{p}{q}\left( y_{\meimej{p}{q},i} \right) \right] \right] \right] \right] \right] \\
						& = \EBL{p}{q}\left( \ovly_{\meimej{p}{q},i} \right) \\
						& = \EBL{p}{q}\left( z_{\meimej{p}{q},i} \right) \in \mfB. \\
					\end{align*}					
		\end{enumerate}
	
	Similarly, by bracketing with elements in $\mfB$, we can realize the indeterminates $z_{\mu,i}^{-1}$ in $Z_f$ as members of the monomials we are considering.  

	Let us return to the task of demonstrating that given a monomial
		\begin{displaymath}
			\Tmanyl,
		\end{displaymath}		
		\begin{itemize}
			\item[]
				$\Evert{i} \left( \Tmanyl \right) \in \mfB$,
			\item[]
				$\Ehort{i} \left( \Tmanyl \right) \in \mfB$,
			\item[]
				$\EUL{i}{j} \left( \Tmanyl \right) \in \mfB$,								
			\item[]
				$\EUR{i}{j} \left( \Tmanyl \right) \in \mfB$, and				
			\item[]
				$\EBL{i}{j} \left( \Tmanyl \right) \in \mfB$.				
		\end{itemize}
	We proceed by induction on $l$.  
	
	For the base case $l = 1$, we know that 
		\begin{displaymath}
			\varphi\left( e_{\mu_1,k_1} \right) = \te_{\mu_1,k_1} 
				= \left\{
					\begin{array}{ll}
						\EUL{p}{q}\left( \ymuk \right) & \text{ if } \mu_1 = \eimej{p}{q} \in \Theta \\
						\EUR{p}{q}\left( \xmuk \right) & \text{ if } \mu_1 = \eipej{p}{q} \in \Omega \\
						\EUR{p}{q}\left( \ymuk \right) & \text{ if } \mu_1 = \eipej{p}{q} \in \Theta \\
						\EBL{p}{q}\left( \xmuk \right) & \text{ if } \mu_1 = \meimej{p}{q} \in \Omega \\
						\EBL{p}{q}\left( \ymuk \right) & \text{ if } \mu_1 = \meimej{p}{q} \in \Theta \\
					\end{array}
				\right.
		\end{displaymath}
	is an element of $\mfB$.  We saw above that after bracketing with appropriate elements in $\mfB$ we get $\EUL{p}{q}\left( \zmuk \right)$ if $\mu_1 = \eimej{p}{q} \in \Theta$, $\EUR{p}{q}\left( \zmuk \right)$ if $\mu_1 = \eipej{p}{q} \in \Theta$, or $\EBL{p}{q}\left( \zmuk \right)$ if $\mu_1 = \meimej{p}{q} \in \Theta$.
	
	So letting $\Tone$ represent $\xmuk$, $\ymuk$, or $\zmuk$, as appropriate, we have $\EUL{p}{q}\left( \Tone \right)$, $\EUR{p}{q}\left( \Tone \right)$, or $\EBL{p}{q}\left( \Tone \right) \in \mfB$ depending on what $\mu_1$ is. 
	
	We also know that 
		\begin{displaymath}
			\varphi\left( f_{\mu_1,k_1} \right) = \tf_{\mu_1,k_1} 
				= \left\{
					\begin{array}{ll}
						\EUL{q}{p}\left( \ymukinv \right) & \text{ if } \mu_1 = \eimej{p}{q} \in \Theta \\
						\EBL{q}{p}\left( \xmukinv \right) & \text{ if } \mu_1 = \eipej{p}{q} \in \Omega \\
						\EBL{q}{p}\left( \ymukinv \right) & \text{ if } \mu_1 = \eipej{p}{q} \in \Theta \\
						\EUR{q}{p}\left( \xmukinv \right) & \text{ if } \mu_1 = \meimej{p}{q} \in \Omega \\
						\EUR{q}{p}\left( \ymukinv \right) & \text{ if } \mu_1 = \meimej{p}{q} \in \Theta \\
					\end{array}
				\right.
		\end{displaymath}
	is an element of $\mfB$.  From the first, third, and fifth of these, we know that bracketing with elements in $\mfB$ leads to, respectively, $\EUL{q}{p}\left( \zmukinv \right)$ if $\mu_1 = \eimej{p}{q} \in \Theta$, $\EBL{q}{p}\left( \zmukinv \right)$ if $\mu_1 = \eipej{p}{q} \in \Theta$, and $\EUR{q}{p}\left( \zmukinv \right)$ if $\mu_1 = \meimej{p}{q} \in \Theta$.  So again, letting $\Tone$ represent $\xmuk$, $\ymuk$, or $\zmuk$, as appropriate, we have $\EUL{q}{p}\left( \Toneinv \right)$, $\EBL{q}{p}\left( \Toneinv \right)$, or $\EUR{q}{p}\left( \Toneinv \right) \in \mfB$, depending on what $\mu_1$ is.  
	
	If $\mu_1 = \eimej{p}{q}$, where $p \neq q$, then
		\begin{displaymath}
			\Evert{p}\left( \Tone \right) = \left[ \EUL{p}{q}\left( \Tone \right), \ \Evert{q}(1) \right] \in \mfB.
		\end{displaymath}
	By Corollary \ref{cor:eia to an}, 
		\begin{displaymath}
			\Evert{i}\left( \Tonem \right) \text{ and } \Ehort{i}\left( \Tonem \right) \in \mfB
		\end{displaymath}
	for every $i = 1, \ldots, r$ and every positive integer $m_1$. Moreover, since for all $j = 1, \ldots, r$, $\Ehort{j}(1)$ and $\Evert{j}(1) \in \mfB$, by Lemma \ref{lemma:off centre quads}, 
		\begin{displaymath}
			\EUL{i}{j}\left( \Tonem \right) \in \mfB,
		\end{displaymath}
	for all $1 \leq i,j \leq r$, with $i \neq j$, and any positive integer $m_1$.  Also
	 	\begin{displaymath}
			\EUR{i}{j}\left( \Tonem \right), \text{ and } \EBL{i}{j}\left( \Tonem \right) \in \mfB,
	 	\end{displaymath}
	 for all $1 \leq i,j \leq r$ and any positive integer $m_1$.  
	 
	We also have 
		\begin{displaymath}
			\Evert{q}\left( \Toneinv \right) = \left[ \EUL{q}{p}\left( \Toneinv \right), \ \Evert{p}(1) \right] \in \mfB.
		\end{displaymath}
	By Corollary \ref{cor:eia to an}, 
		\begin{displaymath}
			\Evert{i}\left( \Tonem \right) \text{ and } \Ehort{i}\left( \Tonem \right) \in \mfB
		\end{displaymath}
	for every $i = 1, \ldots, r$ and every negative integer $m_1$. Since $\Ehort{j}(1)$ and $\Evert{j}(1) \in \mfB$ for all $j = 1, \ldots, r$, we get that
		\begin{displaymath}
			\EUL{i}{j}\left( \Tonem \right) \in \mfB,
		\end{displaymath}
	for all $1 \leq i,j \leq r$, with $i \neq j$, and any negative integer $m_1$, and
	 	\begin{displaymath}
			\EUR{i}{j}\left( \Tonem \right), \text{ and } \EBL{i}{j}\left( \Tonem \right) \in \mfB,
	 	\end{displaymath}
	 for all $1 \leq i,j \leq r$ and any negative integer $m_1$.  	 
		
		If $\mu_1 = \eipej{p}{q}$, where $p \neq q$, then 
			\begin{displaymath}
				\Evert{p}\left( \Tone \right) = \left[ \Ehort{q}(1), \ \EUR{p}{q}\left( \Tone \right) \right] \in \mfB,
			\end{displaymath}
		and
			\begin{displaymath}
				\Ehort{p}\left( \Toneinv \right) = \left[ \EBL{q}{p}\left( \Toneinv \right), \ \Evert{q}(1) \right] \in \mfB.
			\end{displaymath}
		It follows, by the above lemmas and corollaries, that for all $m_1 \in \bbZ$ and for any $1 \leq i, j \leq r$,
			\begin{displaymath}
				\Evert{i}\left( \Tonem \right), \ \Ehort{i}\left( \Tonem \right), \ \EUL{i}{j}\left( \Tonem \right), \ \EUR{i}{j}\left( \Tonem \right), \text{ and } \EBL{i}{j}\left( \Tonem \right) \in \mfB.
			\end{displaymath}

		If $\mu_1 = \meimej{p}{q}$, where $p \neq q$, then 
			\begin{displaymath}
				\Ehort{q}\left( \Tone \right) = \left[ \EBL{p}{q}\left( \Tone \right), \ \Evert{p}(1) \right] \in \mfB,
			\end{displaymath}
		and
			\begin{displaymath}
				\Evert{q}\left( \Toneinv \right) = \left[ \Ehort{p}(1), \ \EUR{q}{p}\left( \Toneinv \right) \right] \in \mfB.
			\end{displaymath}
		Again the above lemmas and corollaries lead to, for all $m_1 \in \bbZ$ and $1 \leq i, j \leq r$,
			\begin{displaymath}
				\Evert{i}\left( \Tonem \right), \ \Ehort{i}\left( \Tonem \right), \ \EUL{i}{j}\left( \Tonem \right), \ \EUR{i}{j}\left( \Tonem \right), \text{ and } \EBL{i}{j}\left( \Tonem \right) \in \mfB.
			\end{displaymath}

	Thus the base case holds.  Next, suppose that for $l \geq 1$ and monomial 
		\begin{displaymath}
			\Tmanyl
		\end{displaymath}
	in $\mfb$, we have, for all $1 \leq i,j \leq r$, 
		\begin{align*}
			& \Evert{i} \left( \Tmanyl \right), \ \Ehort{i} \left( \Tmanyl \right), \ \EUL{i}{j} \left( \Tmanyl \right), \\		
			& \EUR{i}{j} \left( \Tmanyl \right), \text{ and	} \EBL{i}{j} \left( \Tmanyl \right) \in \mfB.				
		\end{align*}
	Having worked through the $l =1$ case, we know that for the monomial $\Tlpom$, $\Ehort{j}\left( \Tlpom \right) \in \mfB$ for all $j = 1, \ldots, r$.  Since $\Evert{i}\left( \Tmanyl \right) \in \mfB$ by the induction hypothesis, we get, by Lemma \ref{lemma:off centre quads}, that
		\begin{displaymath}
			\EUL{i}{j}\left( \Tmanylpo \right) \in \mfB
		\end{displaymath}
	for all $1 \leq i,j \leq r$ with $i \neq j$.  Similarly, by using the induction hypotheses and Lemma \ref{lemma:off centre quads}, we get, for all $1 \leq i,j \leq r$, 
		\begin{displaymath}
			\EUR{i}{j}\left( \Tmanylpo \right), \ \EBL{i}{j}\left( \Tmanylpo \right) \in \mfB.
		\end{displaymath}
	
	By the induction hypothesis, for any $1 \leq i \leq r$, $\Evert{i}\left( \Tmanyl \right) \in \mfB$.  Hence, by Lemma \ref{lemma:eia to hia and meia}, $\EUL{i}{i}\left( \Tmanyl \right) \in \mfB$.  We also know, by the $l=1$ analysis, that, for any $1 \leq i \leq r$, $\Evert{i}\left( \Tlpom \right) \in \mfB$.  But then
		\begin{displaymath}
			\Evert{i}\left( \Tmanylpo \right) = \left[ \EUL{i}{i}\left( \Tmanyl \right), \ \Evert{i}\left( \Tlpom \right ) \right] \in \mfB.
		\end{displaymath}
	
	Since, by the $l=1$ case,  $\Evert{i}\left( \Tlpom \right) \in \mfB$ for any $1 \leq i \leq r$, by Lemma \ref{lemma:eia to hia and meia}, $\EUL{i}{i}\left( \Tlpom \right) \in \mfB$.  We also know, by the induction hypothesis, that for any $1 \leq i \leq r$, $\Ehort{i}\left( \Tmanyl \right) \in \mfB$. But then
		\begin{displaymath}
			\Ehort{i}\left( \Tmanylpo \right) = \left[ \Ehort{i}\left( \Tmanyl \right), \ \EUL{i}{i}\left( \Tlpom \right ) \right] \in \mfB.
		\end{displaymath}
	
	So the result holds for all $l \geq 1$.  
	
	Hence, we have shown that, for all $1 \leq i, j \leq r$, $\sotrpobecc_{-\epsilon_i}$, $\sotrpobecc_{\epsilon_i}$, $\sotrpobecc_{\epsilon_i + \epsilon_j}$, $\sotrpobecc_{-\epsilon_i - \epsilon_j}$, and $\sotrpobecc_{\epsilon_i - \epsilon_j}$ ($i \neq j$) are all contained in $\mfB$.  Moreover, since 
		\begin{displaymath}
			\sotrpobecc_0 = \sum_{\gamma \in \Delta} \left[ \sotrpobecc_{-\gamma}, \ \sotrpobecc_{\gamma} \right],
		\end{displaymath}
	$\sotrpobecc_0$ is also contained in $\mfB$.  We have thus shown that the following proposition holds.	
	\begin{propn}
		The map $\varphi: \gimAd \to \sotrpobecc$ is a surjective Lie algebra homomorphism.
	\end{propn}
	

\subsection{}\label{subsect:is graded}
	
	In this subsection we show that $\varphi: \gimAd \to \sotrpobecc$ is a graded homomorphism and that it induces a map from $\imAd$ to $\sotrpobecc$.
	
	We saw in Sections \ref{sect:im algebras} and \ref{sect:recognition theorem}, respectively, that $\gimAd$ and $\sotrpobecc$ are both $\Gamma$-graded Lie algebras, where
		\begin{displaymath}
			\Gamma = \bigoplus_{\mu \in \Delta} \bbZ \alpha_{\mu}.
		\end{displaymath}
	The map $\varphi: \gimAd \to \sotrpobecc$ is engineered so that, for all $\alpha \in \Gamma$, 
		\begin{displaymath}
			\varphi\left( \gimAd_{\alpha} \right) \subset \sotrpobecc_{\alpha}.
		\end{displaymath}
	That is, the following result holds by design.
	\begin{propn}
		The map $\varphi: \gimAd \to \sotrpobecc$ is also a graded homomorphism.
	\end{propn}
		
	Moreover, since $\sotrpobecc_{\gamma} = 0$ for $\gamma \notin \Delta \cup \{0\}$, we get that the radical $\mfr$ of $\gimAd$ lies in the kernel of $\varphi$.  
	
	\begin{propn}
		There exists a surjective, graded Lie algebra homomorphism 
			\begin{displaymath}
				\phi: \imAd \to \sotrpobecc
			\end{displaymath}
		given by $\phi(u + \mfr) = \varphi(u)$ for any $u + \mfr \in \imAd$, where $u \in \gimAd$.
	\end{propn}


\subsection{}\label{subsect:is central}

	Our work in $\S$ \ref{sect:minimal understanding} revealed that the associative algebra $\mfa$ in the $BC$-graded Lie algebra $\sotrpoaecc$, arising from [ABnG]'s Recognition Theorem, is the algebra $\mfb$ modulo some (possibly zero) ideal $I$ of $\mfb$.  Here $I$, if nonzero, would consist of more relations on the elements of $\mfa$ than we presently have among the elements of $\mfb$.  That is,
		\begin{displaymath}
			\mfa = \mfb / I.
		\end{displaymath}
	In particular, we have a surjective associative algebra homomorphism from $\mfb$ to $\mfa$ that respects the involution and relations.  As a consequence, we get a surjective Lie algebra homomorphism
		\begin{displaymath}
			\sigma: \sotrpobecc \to \sotrpoaecc
		\end{displaymath}
	such that
		\begin{displaymath}
			\sigma \phi = \psi,
		\end{displaymath}
	where
		\begin{displaymath}
			\psi: \imAd \to \sotnpoaecc
		\end{displaymath}
	is the central, surjective, graded Lie algebra homomorphism we established in Section \ref{sect:recognition theorem}.

	But then
		\begin{displaymath}
			\ker \phi \subset \ker \psi \subset \mfz\!\left( \imAd \right),
		\end{displaymath}
	where $\mfz\!\left( \imAd \right)$ denotes the centre of $\imAd$.  The second inclusion, $\ker \psi \subset \mfz\!\left( \imAd \right)$, holds because $\psi$ is a central map.  Thus $\ker \phi \subset \mfz\!\left( \imAd \right)$, implying the following result:
	\begin{propn}
		The map $	\phi: \imAd \to \sotrpobecc$ is a central homomorphism.
	\end{propn}


\section{Future work}

There are at least two avenues of investigation that directly follow from this paper:
		\begin{enumerate}
			\item
				We are studying what $\mfa$, $\eta$, $C$, and $\chi$ look like if we adjoin roots of the form
					\begin{enumerate}
						\item
							$\pm \epsilon_i \in \Delta$, $1 \leq i \leq r$, to a base $\Pi$ of $\Delta_{B_r}$; or							
						\item
							$\pm 2\epsilon_i \in \Delta$, $1 \leq i \leq r$, to a base $\Pi$ of $\Delta_{B_r}$.
					\end{enumerate}
				We are also studying these four components of $\sotrpoaecc$ in the setting where we adjoin a mixture of these short, long, and ``extra-long'' roots in $\Delta$ to a base $\Pi$ of $\Delta_{B_r}$.  Early work suggests that the module $C$ and the hermitian form $\chi$ are no longer trivial.

			\item
				We are looking for realizations of intersection matrix algebras of type $BC_r$, where $r \geq 4$, with a grading subalgebra of type			
					\begin{enumerate}
						\item
							$C_r$, and
						\item
							$D_r$.
					\end{enumerate}
				Here, for example, we expect that $\imAd$ will be a realization of a more general version of $\text{sp}_{2r}(\bbC)$ or $\text{so}_{2r}(\bbC)$, respectively.
		\end{enumerate}

\noindent \emph{Acknowledgements:} The authors would like to thank Frank Lemire and Michael Lau for their helpful discussions.

\vskip20pt

\begin{tabular}{l}
Sandeep Bhargava\\
Dept. of Math. $\&$ Stats. \\
University of Windsor\\
Windsor, Ontario\\
N9B 3P4\\
bhargava@uwindsor.ca
\end{tabular}

\vskip10pt

\begin{tabular}{l}
Yun Gao \\
Dept. of Math. $\&$ Stats. \\
York University \\
4700 Keele Street \\
Toronto, Ontario\\
M3J 1P3 \\
ygao@yorku.ca
\end{tabular}

\end{document}